\documentclass[10pt,a4paper]{article}
\usepackage[latin1]{inputenc} 
\usepackage[english]{babel}
\usepackage[T1]{fontenc} 
\usepackage[all]{xy}
\usepackage{fancyhdr}
\usepackage[dvips]{geometry}
\usepackage{indentfirst,amsfonts,amsmath,amsbsy,amsthm,amssymb,amscd}
\usepackage{enumerate}
\usepackage[normalem]{ulem}  
\usepackage{makeidx}  

\usepackage{color}
\usepackage{mathrsfs}
\usepackage{dsfont}

\usepackage{hyperref}

\usepackage{makeidx}
\usepackage{calligra}
\usepackage[affil-it]{authblk}
\usepackage{graphicx}

\catcode`@=11 \@addtoreset{equation}{section} \catcode`@=12

\newtheorem{theorem}{Theorem}[section]
\newtheorem{lemma}[theorem]{Lemma}
\newtheorem{remark}[theorem]{Remark}
\newtheorem{corollary}[theorem]{Corollary}
\newtheorem{proposition}[theorem]{Proposition}
\newtheorem{definition}[theorem]{Definition}
\newtheorem{example}[theorem]{Example}

\def\XXint#1#2#3{{\setbox0=\hbox{$#1{#2#3}{\int}$ }
\vcenter{\hbox{$#2#3$ }}\kern-.6\wd0}}

\newcommand{\supp}{\mathop{\rm supp}}

\newcommand{\R}{\mathbb{R}}

\newcommand{\eps}{\varepsilon}

\usepackage{epstopdf}

\begin{document}

\title{Nonlinear nonlocal reaction-diffusion problem with local reaction}

\author{An\'{i}bal Rodr\'{i}guez-Bernal\thanks{Partially supported by
    the projects Projects MTM2016-75465 and PID2019-103860GB-I00,
    MINECO, 
Spain and Grupo CADEDIF 
GR58/08, Grupo 920894. email: arober@mat.ucm.es}}
\affil{Departamento de An\'alisis Matem\'atico y Matem\'atica Aplicada\\Universidad Complutense de Madrid, 28040, 
Madrid, Spain \\and \\
Instituto de Ciencias Matem\'aticas\thanks{Partially supported by ICMAT Severo Ochoa
  project SEV-2015-0554 (MINECO)}, CSIC-UAM-UC3M-UCM, Madrid, Spain}
\author{Silvia Sastre-Gomez\thanks{Partially supported by projects MTM2016-75465, MTM2017-83391 and Grupo CADEDIF GR58/08, Grupo 920894.}}
\affil{Departamento de Ecuaciones Diferenciales y An\'alisis Num\'erico\\Universidad de Sevilla, 41012, Sevilla, Spain}

\date{ \today}
%
%

\maketitle

\section{Introduction}

Diffusion is an ubiquitous phenomena in nature. It appears for
example in  the process by which matter is transported from one 
location  of a system to another as a result of random molecular motions, cf. 
\cite{crank_diffusion}.   When the media is reactive to the diffusion process,
reaction diffusion model appear naturally. Hence reaction--diffusion
equations model  different phenomena and   appear in many different
areas such  as physics, biology, chemistry 
and even economics. In biology  for example   they describe 
the evolution in time and space of the density of population of one
or several biological species, \cite{Murray02}. 

When the media is smooth, e.g. a smooth open domain in euclidean space
or a smooth manifold, diffusion is naturally modeled using
differential operators, e.g. Laplace or Laplace--Beltrami operators,
respectively. In such a situation one typicaly encounters the local
reaction--difusion model 
\begin{equation} \label{eq:local_RD}
  u_{t} - \Delta u = f(x,u)
\end{equation}
in some smooth open domain $\Omega \subset \R^{N}$, complemented with some
boundary conditions on the boundary $\partial \Omega$. In
(\ref{eq:local_RD}) the nonlinear term 
\begin{equation} \label{eq:reaction_term}
  f : \Omega \times \R \to \R 
\end{equation}
describes the local rate of production/consumption of the magnitud $u$
at each  point $x\in \Omega$.

In nonsmooth media however diffusion must be described by other means,
and in this context some nonlocal diffusion operators appear
naturally; e.g. \cite{Havlin,fife,hutson, Kigami, Strichartz}. This
approach is applicable in metric measure spaces defined as follows,   \cite{rudin}, 
 \begin{definition}
	A {\bf metric measure space} $\Omega$ is a metric space 
	 $( \Omega,d)$  with a $\sigma$-finite, regular, and complete Borel 
	 measure $dx$ in $ \Omega $, and  that  associates a finite  positive 
	 measure to the balls of $ \Omega $.
\end{definition} 

Then,   if $\Omega$ is a measure metric  space, assume $u(x,t)$ is 
the density of some population at the point $x\in\Omega$ at time $t$, and 
$J(x,y)$ is a positive function defined in $\Omega\times\Omega$ that 
represents the fraction of the population  jumping from a location $y$
to location  $x$, per unit time. Then  $\int_{\Omega} J(x,y)u(y,t)\, dy$ is the rate at which the individuals arrive to 
location $x$ from all other locations $y\in\Omega$. Analogously,
$\int_{\Omega} J(y,x)u(x,t)\, dy = u(x,t) \int_{\Omega} J(y,x)\, dy$ is the rate at which individuals
leave from location $x$ to any other place in $\Omega$. 
%
%
Hence  the evolution in time of the population 
can be written as 
\begin{equation}\label{intro_linear_diffusion}
	\left\{
	\begin{array}{rl}
	u_t(x,t) \!\!\!\!& \displaystyle =\int_{\Omega}J(x,y)u(y,t)\,
                           dy\!- h_{*}(x)u(x,t),\quad x\in\Omega,\, t>0,\\
	u(x,0) \!\!\!\! & =u_0(x), \quad x\in\Omega,
	\end{array}
	\right.
\end{equation}
where $u_0$ is the initial distribution density of the population and $h_{*}(x) =
\int_{\Omega} J(y,x)\, dy$. Observe that in a symmetric media we have
$J(x,y)=J(y,x)$ and (\ref{intro_linear_diffusion}) can be recast as 
\begin{displaymath}
 	u_t(x,t) =\int_{\Omega}J(x,y) \Big( u(y,t) - u(x,t) \Big) \,
                           dy ,\quad x\in\Omega,\, t>0 . 
\end{displaymath}

This allows  to study 
diffusion processes in very different types of media  like, for example, graphs 
(which are used to model complicated structures in chemistry, molecular 
biology or electronics, or they can also represent basic electric
circuits in digital computers), compact manifolds, multi-structures composed 
by several compact sets with different dimensions (for example, a dumbbell 
domain), or even some fractal sets such as the Sierpinski gasket, see
\cite{ARB_SSG_16} for some details. The case when $\Omega$ is an open
set of euclidean space (\ref{intro_linear_diffusion}) and variations of it  have been consider thoroughly in 
\cite{Rossi_libro, BatesChmaj, Bates_97,  Coville_2016, ref4, Coville_2006, Coville_2010, Coville_2008, Coville_2005, Shen_2015, Shen_2010, Shen_2012} and references therein. 
Other approaches to diffusion in nonsmooth media can be found in 
\cite{Havlin, Lapidus, Kigami, Strichartz}. 

A nonlinear version of (\ref{intro_linear_diffusion}) that we consider
in this paper reads 
\begin{equation}\label{intro_nonlinear_diffusion}
	\left\{
	\begin{array}{rl}
	u_t(x,t) \!\!\!\!& \displaystyle =\int_{\Omega}J(x,y)u(y,t)\,
                           dy\! +  f(x,u(x,t)) ,\quad x\in\Omega,\, t>0,\\
	u(x,0) \!\!\!\! & =u_0(x), \quad x\in\Omega,
	\end{array}
	\right.
\end{equation}
where   the nonlinear local  reaction is as in
(\ref{eq:reaction_term}).


The aim of this paper is to make a general 
study of the nonlinear problem (\ref{intro_nonlinear_diffusion}) and show some similarities 
and differences between (\ref{intro_nonlinear_diffusion}) and the local reaction-diffusion 
problem (\ref{eq:local_RD}). We will show that both models share positivity properties and the 
strong maximum principle. However, the lack of smoothing effect  of the
 linear equation (\ref{intro_linear_diffusion})   affects to some results of existence and asymptotic behaviour of 
the nonlocal problem which are weaker than for the local problem. 


The paper is organized as follows. In
section \ref{sec:prelim-line-equat} we introduce all the standing
assumptions on the metric measure space $\Omega$ and the nonlocal
kernel $J$. Also we recall some of the results in
\cite{RB2017_max_pples, ARB_SSG_16} on  the linear nonlocal
problem (\ref{intro_linear_diffusion}) that will be used throughout
this paper. In particular, depending on properties of $J$ we can
consider initial data for (\ref{intro_linear_diffusion}) and
(\ref{intro_nonlinear_diffusion}) in the spaces 
$X=L^p(\Omega)$, $1\le p\le \infty$, or $X=\mathcal{C}_b(\Omega)$. 

In Section \ref{exist_uniq_pos_nonlinear} we prove global existence of
solutions for (\ref{intro_nonlinear_diffusion}) for initial  $u_{0} \in X$ for some classes of
nonlinear terms $f(x,u)$. First, in Section 
\ref{exist_uniq_pos_Glob_Lipschitz}  we consider the case when $f$ is
globally Lipschitz and  $u_{0} \in X$. 
%
%
We also prove that solutions of (\ref{intro_nonlinear_diffusion})
satisfy both weak, strict or strong comparison and maximum principles,
depending on conditions for the kernel $J$. These results are analogous to the results of  the local 
nonlinear reaction-diffusion problem with boundary conditions, (cf.  
\cite{arrieta_carvalho_anibal}).   
Using these results, in  Section \ref{existence_sol_with_f_loc_Lip} we
first prove global existence,  
uniqueness, comparison and maximum principles for  the solution of
(\ref{intro_nonlinear_diffusion}) for bounded initial data  when the
nonlinear term $f(x,u)$, is locally Lipschitz in the variable  
$u\in\mathbb{R}$, uniformly with respect to $x\in\Omega$,  and
satisfies a 
 sign condition that reads 
 \begin{equation}\label{intro_structure_condition_f}
	f(x,u)u\le Cu^2+D|u|,\quad \mbox{ for all } \ u \in \R,  \ x\in\Omega.
      \end{equation}
for some  $C,D\in\R$ with $D>0$. 
      This differs from the local reaction-diffusion problems since for 
$f$ locally Lipschitz the local existence, at least for smooth
initial data,  can be proved without any  extra sign condition on
$f$. Finally, for some particular nonlinear terms that satisfy
(\ref{intro_structure_condition_f}) and some growth condition we are
able to prove global existence for initial data $u_{0}$ in a suitable
$L^{p}(\Omega)$ space. In this  case, we also  prove comparison and maximum
principles for the solutions. 


In Section \ref{asymptotic_estimates}, we give some asymptotic
estimates of the solution constructed in Section
\ref{exist_uniq_pos_nonlinear}. In particular we prove asymptotic
pointwise $L^{\infty}(\Omega)$ bounds on the solutions. These
estimates are improved in  Section \ref{section_extr_equil} where we prove
the existence of two extremal ordered  equilibria $\varphi_m \leq \varphi_M$ in
$L^{\infty}(\Omega)$ of (\ref{intro_nonlinear_diffusion}).  These
extremal equilibria give pointwise asymptotic   bounds of any weak limit in $L^p(\Omega)$ for
$1\le p<\infty$, or weak* limit in $L^{\infty}(\Omega)$, of the
solution of (\ref{intro_nonlinear_diffusion}) with initial data
$u_0\in L^p(\Omega)$ or $ L^\infty(\Omega)$, respectively. We also prove
that the maximal extremal equilibria is ``stable from above'' and the
minimal extremal equilibria is ``stable from below''.  
We find again here  another difference with the nonlinear
local problem (\ref{eq:local_RD}), where the asymptotic dynamics of
the solution enter between  extremal equilibria, uniformly in
space, for bounded sets of initial data, cf.
\cite{Anibal_Alejandro_extremal_equilibria} and they are part of the
global attractor of the problem,. This difference is again  due to
the lack of smoothing  of the linear group associated to
(\ref{intro_linear_diffusion}).
Another striking difference with
(\ref{eq:local_RD}) is given in Example
\ref{example_discontinuous_equilibria} where we show that it is
possible to construct an uncountable family of nonisolated, discontinuous,
piecewise constant equilibria that even may coincide in sets of
positive measure. Such family of equilibria can be made also of
positive equilibria. Observe that the discontinuity of equilibria is
related once more to the lack of smoothing of solutions of
(\ref{intro_nonlinear_diffusion}). 

Then in  Section \ref{sec:nonnegative-solutions} for   nonnegative solutions we prove that
if $f(x,0)$ is nonnegative then there exists a minimal nonnegative
equilibrium, $0\le \varphi_m^+\le \varphi_M$, such that the solutions
of (\ref{intro_nonlinear_diffusion}) with nonnegative initial data
enter asymptotically between these nonnegative extremal equilibria. If
moreover $f(x,0)=0$ and $u=0$ is linearly unstable then every
nontrivial nonnegative equilibria is strictly positive. We finally
give a sufficient conditions for uniqueness of a positive
equilibrium. In such a case we obtain that the unique equilibria is
globally asymptotically stable for nonnegative solutions.
Then  in Section \ref{logistic_nonlin} we  also  
analyze in detail the case of logistic nonlinearities for which
\begin{displaymath}
  f(x,u) = g(x) + n(x) u - m(x) |u|^{\rho-1} u
\end{displaymath}
with $g, n, m\in L^{\infty}(\Omega)$, $m\ge 0$ not identically zero
and $\rho>1$. For these problems we show that extremal equilibria
always exist and give conditions on $g,m,n$ that guarantee uniqueness
of a positive, globally stable, positive equilibria. In doing so we
prove a result that  states that  by acting on an arbitrary small
subset of the domain with a large negative constant,  we can shift the
spectrum of a linear nonlocal operator plus a potential $h(x)$, 
to have negative real part, a result that does not hold for the local
diffusion operator $-\Delta + h(x)I$.

Finally, in Section \ref{asymp_beh} we give some further comments
about the asymptotic behaviour of
\eqref{intro_nonlinear_diffusion} and the lack of asymptotic
compactness. In particular, we show that if we had enough compactness
to guarantee that, along subsequences of  large times, solutions converge a.e. $x\in
\Omega$, then it would be possible to prove the solutions approach
equilibria and even the existence of a suitable attractor.  We have
failed in proving that such pointwise asymptotic convergence holds
true for (\ref{intro_nonlinear_diffusion}).

%

\section{Preliminaries on linear equations} 
\label{sec:prelim-line-equat}

Let $\Omega$ be a metric measure space 
and let $J$ be a nonnegative  kernel  defined as $J:\Omega\times
\Omega\to\R$, considered as a mapping 
\begin{displaymath}
\Omega\ni x\mapsto J(x,\cdot)\ge 0 .
\end{displaymath}
Then  consider the  nonlocal diffusion operator  given by  
\begin{displaymath}
  Ku(x)=\int_{\Omega}J(x,y)u(y)\, dy,  \quad x\in \Omega
\end{displaymath}
for  suitable functions defined in $\Omega$. 
Consider also $h \in L^{\infty}(\Omega)$ a bounded measurable
function. 

The operator $K$ will be considered below in the Lebesgue spaces
$X=L^p(\Omega)$, with $1\le p\le\infty$ or in the space of bounded
continuous functions $X= C_b(\Omega)$. In the latter case  we
will assume $h\in C_b(\Omega)$.

As general notations, $|A|$ will denote the measure of a measurable
set $A\subset \Omega$.

%
%
%
%
%
%
%

\subsection{Stationary  problems}

The next results state regularity and compactness  properties of $K$ 
derived from properties of the kernel $J$, see e.g. \cite{ARB_SSG_16}
for details. First we have the following.

\begin{proposition} \label{prop:KJ_in_LXX}

\noindent i) 
Let $1\le p\le \infty$, if $J\in	L^p(\Omega, L^{p'}(\Omega))$ then  
	$K\in \mathcal{L}(L^p(\Omega), L^p(\Omega))$. 
If moreover, $p<\infty$ then $K\in \mathcal{L}(L^p(\Omega),
L^p(\Omega))$ is compact. 

\noindent ii)    If $J\in	 C_b(\Omega, L^1(\Omega))$ then  
	$K\in \mathcal{L}( C_b(\Omega),  C_b(\Omega))$. 
Moreover if $J\in BUC(\Omega,  L^{1}(\Omega))$, then $K\in \mathcal{L}(L^\infty(\Omega), 
 C_b(\Omega))$ is compact. In particular,  $K\!\!\in \!\mathcal
{L}(L^\infty(\Omega), L^{\infty}(\Omega))$ is compact and $K\!\!\in \!
\mathcal{L}( C_b(\Omega),  C_b(\Omega))$ is compact.

\noindent iii)  
If $|\Omega|<\infty$ and
$J\in L^{\infty}(\Omega, L^{p_{0}'} (\Omega))$ for some
$1\leq p_{0} \leq \infty$, then
$K\in \mathcal{L}(L^p(\Omega), L^p(\Omega))$, for all
$p_{0} \le p < \infty$ and is compact.
If moreover, $J\in BUC(\Omega, L^{p_{0}'} (\Omega))$ then $K\in
\mathcal{L}(L^{\infty}(\Omega), L^{\infty}(\Omega))$ and $K\in \mathcal{L}(
C_b(\Omega),  C_b(\Omega))$ are compact.

\end{proposition}

Notice in particular  that if
$J \in L^{\infty}(\Omega, L^{1}(\Omega))$ then $K\in
\mathcal{L}(L^{\infty}(\Omega), L^{\infty}(\Omega))$ and  the function 
\begin{equation} \label{eq:definition_h0}
 h_0(x)=\int_{\Omega}J(x,y)dy   
\end{equation}
satisfies $h_0\in L^{\infty}(\Omega)$. If moreover
$J\in BUC(\Omega,L^1(\Omega))$ then $h_0\in C_{b}(\Omega)$.

Also,  for a measurable function $g:\Omega \to \R$ we define the
{\bf essential range} of $g$ (range for short)  as 
\begin{equation} \label{eq:essential_range}
  R(g) = \{s\in \R: \ |\{x: \ |g(x)-s|<\eps\}|>0 \  \mbox{for all} \ \eps>0\}
\end{equation}
which coincides with the set of $s \in \R$ such that $\frac{1}{g(x)-s}
\in L^{\infty}(\Omega)$. Also, if $g$ is continuous this coincides
with the image set of $g$. 
We will also make use of the essential infimum and essential
supremum of a measurable function, which we will denote infimum and
supremum for short, defined as 
\begin{displaymath}
  \inf_{\Omega} g= \sup\{ \alpha \in \R: \ |\{ g \leq \alpha\}|
  =0\} , \quad 
  \sup_{\Omega} g = \inf\{ \alpha \in \R: \ |\{ g \geq \alpha\}|  =0\} . 
\end{displaymath}
Note that  both $\inf_{\Omega} g$ and $\sup_{\Omega} g$ belong to the
(essential) range $R(g)$ if they are finite.


Then Proposition \ref{prop:KJ_in_LXX} implies that the spectrum of
$K-hI$ satisfies $\sigma(K-hI) = \sigma_{ess} \cup \sigma_{p}$ where
the essential essential spectrum is
\begin{displaymath}
  \sigma_{ess}= R(-h) 
\end{displaymath}
where $R(-h)$ is the essential range of the function $-h$, see
(\ref{eq:essential_range}), and a (possibly empty) discrete  point spectrum
$\sigma_{p}  = \left\{\mu_n\right\}_{n=1}^{M}$,
$M\in\mathbb{N}\cup\{\infty\}$. If $M=\infty$, then
$\; \left\{\mu_n\right\}_{n=1}^{\infty}\;$ accumulates in
$R(-h)$.



Note that the essential spectrum $\sigma_{ess}(K-hI)= {R(-h)}$ is
independent of $X$ and it is formed by the points such 
that $K-(h+\lambda)I$ is not a
Fredholm operator of index zero. Also note that the point spectrum
$\sigma_{p}(K-hI) = \left\{\mu_n\right\}_{n=1}^{M}$ is potentially
dependent of the space $X$. Hence, the following result, taken from
\cite[Proposition 3.25]{ARB_SSG_16}, guarantees that the point
spectrum, hence the whole spectrum $\sigma(K-hI)$,  is  independent of $X$. 

\begin{proposition}\label{independent_spectrum_K_h}

 Assume $|\Omega|<\infty$ and
  $J\in L^{\infty}(\Omega, L^{p_{0}'} (\Omega))$ for some $1\leq p_{0}
  \leq \infty$ and
  $h\in L^{\infty}(\Omega)$ then for all $p_{0} \leq p< \infty$,
  $K-hI\in\mathcal {L}(L^{p}(\Omega),L^{p}(\Omega))$, and
  $\sigma_{L^p(\Omega)} (K-hI)$ is independent of $p$.

  If moreover $J \in BUC(\Omega, L^{p_{0}'}(\Omega))$, the spectrum
  above coincides also with $\sigma_{L^{\infty}(\Omega)} (K-hI)$. If,
  additionally,  $h\in C_{b}(\Omega)$, the spectrum above coincides also
  with $\sigma_{C_{b}(\Omega)} (K-hI)$.

\end{proposition}

  Below we give several definitions that will be useful for the following results. 
  
 \begin{definition}
	Let $z$ be a  nonnegative 	measurable  function $z:\Omega\to\R $. We define the  
	{\bf essential support} of  $z$ (support for short) as: 
	\begin{displaymath} 
          \supp(z) =\big\{x\in\Omega\, : \;\forall \delta>0,\;
          |\{y\in\Omega:\;z(y)>0\}\cap B(x,\delta)|>0\},
	\end{displaymath}
 	where $B(x, \delta)$ is the ball 
	centered in $x$, with radius $\delta$.
\end{definition}
%
%
Observe that for a measurable nonnegative  function $z:\Omega\to\R $ 
\begin{displaymath}
  \supp(z)=\Omega \ \mbox{if and only if} \ z>0 \ \mbox{a.e. in
    $\Omega$}.  
\end{displaymath}
In such a case we say that $z$ is
{\bf essentially positive} (positive for short), and write it $z>0$.

Given two measurable functions $w, z:\Omega\to\R $ we
will say that $w$ is (essentially) strictly above  $z$ and write
$w>z$, if $w-z>0$ in the sense above.

We also define 
\begin{definition}\label{R_connected}
If $R>0$, we 	say that $\Omega$ is 
	{\bf $R$-connected} if for all 
	$x,y\in\Omega$, there exist $ N\in \mathbb{N}$
	and a finite set of points $\{x_0,\dots,x_{N}\}$ in $\Omega$ such that 
	$x_0=x$, $x_N=y$ and 
	$d(x_{i-1},x_{i})<R$, for all $i=1,\dots,N$.
\end{definition}

\begin{definition}
Assume that $J$, $h$ and $X$  are as in Proposition
\ref{independent_spectrum_K_h}. 

\noindent i)   We say $\mu \in \R$ is a {\bf principal eigenvalue} of $K-hI$ in $X$ iff
  there exists $0<\phi \in X$ such that 
  \begin{displaymath}
    K\phi -h \phi = \mu \phi  \quad  \text{in} \ \Omega  . 
  \end{displaymath}

\noindent ii) We say that for $\lambda \in \R$ the {\bf maximum principle}
is satisfied if $u\in X$ with 
\begin{displaymath}
  K u - (h+\lambda) u \leq 0\ \text{in} \ \Omega , \quad
  \text{implies} \quad u \geq 0 \
  \text{in} \ \Omega . 
\end{displaymath}

\noindent iii) We say that for $\lambda \in \R$ the {\bf strong maximum principle}
is satisfied if $u\in X$ with 
\begin{displaymath}
  K u - (h+\lambda) u \leq 0\ \text{in} \ \Omega , \quad
  \text{implies either $u=0$ or $u \geq \alpha >0$} 
\end{displaymath}
for some $\alpha >0$. 

\end{definition}

The following theorem gives sufficient conditions for the existence of
the principal eigenvalue of $K-hI$, and it gives a characterization of
the principal eigenvalue when the measure of $\Omega$ is finite, cf. \cite{RB2017_max_pples}. 

\begin{theorem} \label{thr:Lambda}

  Assume that $J$, $h$ and $X$ are as in Proposition \ref{independent_spectrum_K_h},
  $\Omega$ is $R$--connected, $|\Omega| < \infty$,
\begin{equation} \label{eq:measure_non_degenerates}
  |B(x,R)| \geq \mu_{0} >0 \quad \text{for all} \quad x\in \Omega
  , 
\end{equation}
$J\geq 0$ and  there exists $J_{0} >0$ 
\begin{equation}\label{eq:J_strict_positive}
  J(x,y) >  J_{0} >0 \quad \mbox{ for all }x,\,y\in\Omega\mbox{, such
    that \  }d(x,y) < R 
\end{equation}
%
%
and $J \in L^{\infty}(\Omega,
L^{p'}(\Omega))$ so  $K\in \mathcal{L} (X,
L^{\infty}(\Omega))$.

Then $   \Lambda = \sup Re(\sigma_{X}(K-hI)) $ 
can be characterized as 
  \begin{equation}\label{eq:alternative_4_Lambda}
    \Lambda = \inf_{0< \varphi \in X} \sup_{\Omega} \frac{K
      \varphi -h \varphi}{\varphi}  = \sup_{0< \phi \in X} \inf_{\Omega} \frac{K
      \phi -h \phi}{\phi} 
  \end{equation}
%
and satisfies $-\inf_{\Omega} h  \leq \Lambda \leq 
\sup_{\Omega} \big(h_{0} -h \big) $ 
where $h_0$  is defined in (\ref{eq:definition_h0}). 
In particular, $\Lambda$ is the only possible principal eigenvalue
of $K-hI$ in $X$.

\noindent i) If $  \Lambda > -\inf_{\Omega} h$ 
then $\Lambda$ is the  principal eigenvalue of $K-hI$ in $X$. 
In such  a case $\Lambda$ is a simple isolated eigenvalue of
$K-hI$ in $X$ with bounded eigenfunctions.
If moreover $J
\in BUC(\Omega, L^{p_{0}'}(\Omega))$ and $h\in C_b(\Omega)$, the
eigenfunction is continuous. 


\noindent ii) The maximum principle is satisfied for $\lambda
>\Lambda$ and is not satisfied for $\lambda < \Lambda$ nor for
$\lambda= \Lambda > -\inf_{\Omega} h$. 

\noindent iii) If $\lambda >\Lambda$ 
then the strong maximum principle is satisfied.

\end{theorem}

Denoting $m=\inf_{\Omega} h$,   some criteria were also  developed in
\cite{RB2017_max_pples}  to guarantee that  $\Lambda >
-m$, hence $\Lambda$ is the principal eigenvalue of $K-hI$. These
include either one of the following conditions 
\begin{equation} \label{eq:measure_h=m}
|\{h = m \}| >0,    
\end{equation}
or, there  exists $x_{0} \in \Omega$ and $r\leq R$ as in
(\ref{eq:J_strict_positive}) such that 
\begin{equation} \label{eq:1_over_h-m_notintegrable}
\int_{B(x_{0}, r)}  \frac{ dx}{h(x) -m} 
= \infty , 
\end{equation}
or, 
\begin{equation} \label{eq:oscilation_h}
  osc_{\Omega}(h):= \sup_{\Omega} h - \inf_{\Omega} h < \inf_{\Omega}
  h_{0}  
\end{equation}
where $h_0$  is defined in (\ref{eq:definition_h0}). 

The following result, 
for the case of a
symmetric kernel, gives an alternative description of $\Lambda$
using the variational properties in $L^{2}(\Omega)$. 

\begin{proposition} \label{prop:spectrum_as_in_L2}
  
Assume $\Omega$, $J$ and $h$ are as in Theorem \ref{thr:Lambda}
and assume furthermore that 
$J\in L^{\infty}(\Omega, L^{p_{0}'} (\Omega))$ for some  $1\leq p_{0}
\leq 2$ and $J(x,y)=J(y,x)$. 

Then the spectrum of $K -hI$ is real, independent of $X$  and 
        \begin{displaymath}
          \Lambda = \sup\limits_{\varphi \in L^2(\Omega) \atop
	  \|\varphi\|_{L^2(\Omega)}=1}  E(\varphi)
        \end{displaymath}
	where 
	\begin{displaymath}
	  E(\varphi)=-\frac{1}{2}\int_{\Omega}\int_{\Omega}J
	(x,y)(\varphi(y)-\varphi(x))^2 \, dy\,dx - \int_{\Omega} \left( h(x)-h_{0}(x)
	\right)
	    \varphi^{2}(x)\, dx
	\end{displaymath}
with  $h_{0}$ as in  (\ref{eq:definition_h0}).  

\end{proposition}

Also in \cite{RB2017_max_pples} the following criteria for the sign of the principal eigenvalue $\Lambda$ were
obtained. Notice that this information will be used for the stability
of the evolution equations, see Proposition
\ref{prop:intro_behavior_stability} below. 

\begin{proposition} \label{prop:citeria_4_sign_Lambda}

With the assumptions in Theorem \ref{thr:Lambda} and denoting
$m=\inf_{\Omega} h$, we have the following results.

\noindent i)  If $m <0$ then $\Lambda >0$.

\noindent ii) If $m=0$, and  either $|\{h=0\}| >0$ or $\frac{1}{h}
\notin L^{1}_{loc}(\Omega)$ or
$\sup_{\Omega} h < \inf_{\Omega} h_{0}$ or $h + \delta \leq h_{0}$
for  some $\delta >0$, then $\Lambda >0$ and is the principal
eigenvalue. 


\noindent iii) If $m >0$, assume  there exists $0<\xi \in X$ such
that $K\xi -h \xi \nleq 0$,   then  $\Lambda <0$.

\noindent iv) If $m >0$, assume  there exists $\eta  \in X$ that
changes sign in $\Omega$ such that $K\eta -h \eta \leq 0$ then $\Lambda >0$.

\end{proposition}

We also get the following result that sets $h_{0}$ as a
threshold for the sign of $\Lambda$. 

\begin{corollary} \label{cor:theshold_h0}
  
\noindent i)   If $h=h_{0}$ then  $\Lambda =0$ with constant eigenfunction.  

\noindent ii)   If $h_{0}\lneqq h$ then $\Lambda <0$.

\noindent iii)  If for some $\delta >0$, $h+\delta \leq h_{0}$ then
$\Lambda>0$.

\noindent iv) In the symmetric case and $1\leq p_{0} \leq 2$ as in
Proposition \ref{prop:spectrum_as_in_L2}, then $\int_{\Omega} h <
\int_{\Omega} h_{0}$ implies  $\Lambda >0$. 

\end{corollary}

\subsection{Evolutionary problems} 

In this section we present results concerning existence, uniqueness,
maximum principles, bounds and stability results concerning linear
evolution  equation \eqref{intro_linear_diffusion}, cf. \cite{RB2017_max_pples, ARB_SSG_16}. For more information see \cite{Rossi_libro, Rossi, Coville_2008}.

\begin{proposition}\label{positive_solution_with_u_0_positive}

Assume that
$J\in L^{p}(\Omega, L^{p'}(\Omega))$  for some  $1\leq p \leq \infty$ and then denote $X=
L^{p}(\Omega)$ and assume $h\in L^{\infty}(\Omega)$. Alternatively,
assume  $J\in C_{b}(\Omega, L^{1}(\Omega))$ and then denote $X= C_b(\Omega)$ and  assume $h\in C_b(\Omega)$.


\noindent i) Then $L=K-hI\in\mathcal{L}(X,X)$  generates a group 
$e^{Lt}\in\mathcal{L}(X,X)$ for $ t\in\mathbb{R}$,
and the solutions of the initial value problem
	 \begin{equation}\label{linear_equation}
	  \left\{
	 \begin{array}{lll}
	 u_t(x,t) & =(K-hI)u(x,t) & x\in\Omega,  \\
	 u(x,0) & =u_0(x) & x\in\Omega, \end{array}
	 \right.
	 \end{equation}
are given by $u(t)= e^{Lt}u_0$. Finally if  $\Lambda= \sup
Re(\sigma_{X}(L)) < \delta$ then $\|e^{Lt}\|_{\mathcal{L}(X)} \leq M
e^{\delta t}$. 

\noindent ii)  For each $u_{0} \in X$ and $t\in \R$ we have 
\[
 u(x,t)=  e^{Lt}u_0 (x)=e^{-h(x)t}u_0(x)+\displaystyle\int_{0}^t e^{-h(x)(t-s)}K(u) (x,s)ds.
 \]

\noindent iii) 	If  $J\ge 0$, then 
	 for every nonnegative  $u_0\in X$, the solution of  problem
         (\ref{linear_equation}),  $u(t)= e^{Lt}u_0$, is nonnegative for  all $t\ge 0$, and it is nontrivial if $u_0\not\equiv 0$. 
	
Moreover, if $J$ satisfies 
\begin{equation}
  \label{eq:J_positive}
  J(x,y)>0\, \mbox{ for all } \,x,\,y\in\Omega,\, \mbox{ such that } d(x,y)\!<\!R,
\end{equation} 
	for some $R>0$
        and $\Omega$ is $R$-connected as in Definition
        \ref{R_connected}, then  for every $u_0\in X$, nonnegative and
        not 	identically zero,  
        \begin{displaymath}
          \supp(e^{Lt}u_0) = \Omega, \ \text{for all $t>0$}, 
        \end{displaymath}
that is, the solution  of (\ref{linear_equation}) is strictly positive in $\Omega$, for all $t>0$. 
 \end{proposition}

In \cite{RB2017_max_pples} the following results  were proven
concercing  the strong maximum principle and some stability results.  


\begin{proposition} {\bf (Parabolic strong maximum principle)} 
\label{prop:Linear_Parab_Strong_MaxPple}

With the assumptions in Proposition
\ref{positive_solution_with_u_0_positive}, assume furthermore that  $|\Omega|<\infty$, $\Omega$ is  $R$--connected and  the measure satisfies
  (\ref{eq:measure_non_degenerates}) and $J$ satisfies
  (\ref{eq:J_strict_positive}). \\
Then for every $u_0\in X$, nonnegative and  not identically zero,
\begin{displaymath}
  \inf_{\Omega} e^{Lt}u_0  >0, \quad t>0 . 
\end{displaymath}

\end{proposition}



	

\begin{proposition} {\bf (Bounds and stability)}
\label{prop:intro_behavior_stability}

  Assume $\Omega$ is  $R$--connected,   $|\Omega| < \infty$,
(\ref{eq:measure_non_degenerates}) holds true, $J$  satisfies
(\ref{eq:J_strict_positive}). Also assume  $J \in L^{\infty}(\Omega,
L^{p'}(\Omega))$ with $1\leq p \leq \infty$   and then denote $X=
L^{p}(\Omega)$ and assume $h\in L^{\infty}(\Omega)$. Alternatively,  assume $J\in
C_{b}(\Omega, L^{1}(\Omega))$ and then denote $X= C_b(\Omega)$ and  assume $h\in C_b(\Omega)$.


Fix any  $  \tilde{\lambda} < \Lambda < \lambda$.
%
Then 

\noindent i) Any solution of (\ref{intro_linear_diffusion}) with $u_{0}\in X$
satisfies 
\begin{displaymath}
  \|u(t)\|_{X} \leq M e^{\lambda t} \|u_{0}\|_{X} , \quad t\geq 0. 
\end{displaymath}

\noindent ii) Assume either $\Lambda >-\inf_{\Omega} h$ or $J\in
BUC(\Omega, L^{p'}(\Omega))$. Also, by Proposition
\ref{prop:Linear_Parab_Strong_MaxPple}, assume without loss of
generality that $0\leq u_{0} \in X$ is such that $u_{0} \geq
\alpha>0$.  

Then there
exists a positive bounded function $\tilde\varphi$ in $\Omega$ such that 
\begin{displaymath}
0 <   e^{\tilde{\lambda} t}  \tilde{\varphi} (x) \leq u(x,t) , \quad
x\in \Omega, \ t>0. 
\end{displaymath}

\noindent iii) For any solution of (\ref{intro_linear_diffusion})  with $u_{0}
\in L^{\infty}(\Omega)$ there exists a positive function $\varphi \in
X$ such that 
\begin{displaymath}
  |u(x,t,u_{0})| \leq  e^{\lambda t} \varphi(x) \quad x\in \Omega, \ t>0. 
\end{displaymath}

Both parts ii) and iii) hold true for $\lambda=\tilde{\lambda}=
\Lambda$ provided $\Lambda > -\inf_{\Omega} h$. 

In particular, if $\Lambda <0$ all solutions of
(\ref{intro_linear_diffusion})  converge to $0$ in $X$ as $t\to 
\infty$. Moreover, if $u_{0}\in L^{\infty}(\Omega)$ then $u(t) \to 0$
uniformly in $\Omega$ as $t \to \infty$.

On the other hand, if  $\Lambda >0$ then all positive solutions of
(\ref{intro_linear_diffusion})  converge pointwise to  $\infty$ as
$t\to \infty$.

\end{proposition}

\begin{remark}\label{rem:lower_bounds}
For later use, notice that part  ii)   is based on the fact  that the spectrum of $K-hI$ in $X$ coincides
with the spectrum in $L^{\infty}(\Omega)$ and then from
(\ref{eq:alternative_4_Lambda}) given $\tilde{\lambda} < \Lambda$ we can
chose $0<\tilde{\varphi} \in L^{\infty}(\Omega)$ such that
$\tilde{\lambda} <  \inf_{\Omega} \frac{K      \tilde{\varphi} -h 
    \tilde{\varphi}}{\tilde{\varphi}}     \leq  \Lambda$ 
and then  $\underline{u}(x,t) = e^{\tilde{\lambda} t} \tilde{\varphi}
>0$ satisfies 
$ \underline{u}_{t} \leq K \underline{u} -h \underline{u}$, 
i.e. it is a  positive subsolution of (\ref{linear_equation}). 
\end{remark}

\section{Existence, uniqueness, positiveness and comparison results 
for nonlinear problems} 
\label{exist_uniq_pos_nonlinear}

In this section we prove results of existence and uniqueness for some 
nonlocal nonlinear problems of the form
\begin{equation}\label{nonlinear_diffusion}
	\left\{
	\begin{array}{rl}
	u_t(x,t) \!\!\!\!& \displaystyle =\int_{\Omega}J(x,y)u(y,t)\,
                        dy\!  - h(x) u(x) +  f(x,u(x,t)) ,\quad x\in\Omega,\, t>0,\\
	u(x,0) \!\!\!\! & =u_0(x), \quad x\in\Omega 
	\end{array}
	\right.
      \end{equation}
      for some classes of locally Lipschitz functions $f: \Omega \times \R \to \R $.
We will prove also maximum principles and comparison results for
solutions with initial data $u_{0} \in X$ with either
$X=L^{p}(\Omega)$, $1\leq p\leq \infty$ or $X=C_{b}(\Omega)$,
depending on $J$ and $h$ as in Proposition \ref{positive_solution_with_u_0_positive}. For such
$f$ we will consider the associated Nemitcky operator, defined for
mesurable functions $u$ defined in $\Omega$ as
\begin{equation}\label{eq:Nemitcky}
  F(u)(x)=f(x,u(x)), \quad x\in \Omega . 
\end{equation}


Now  we prove some monotonicity and comparison properties for the
problem (\ref{nonlinear_diffusion}). For this we define the
following.

\begin{definition} 
For $u_{0} \in X$ denote by $u(t,u_{0}, f)$ the solution of
(\ref{nonlinear_diffusion}) which we assume to exist for some class
of nonlinear terms $f$.

\noindent i) 
We say that (\ref{nonlinear_diffusion}) satisfies a weak
comparison principle if for any   $f_{0}\geq f_{1}$ and 
$u_0,\,u_1\in X$ such that 	$u_0\ge u_1$, then
\begin{displaymath}
u( t, u_{0},f_{0})\ge u(t, u_{1}, f_{1}),\; \mbox{ for all }\;t\ge0,  
\end{displaymath}

\noindent ii) If moreover $u_{0} \neq u_{1}$ or  $f_{0} \neq f_{1}$
and 
\begin{displaymath}
u( t, u_{0},f_{0})> u(t, u_{1}, f_{1}),\; \mbox{ for all }\;t\ge0,  
\end{displaymath}
we say that (\ref{nonlinear_diffusion}) satisfies a strict 
comparison principle.

\noindent iii) If furthermore 
\begin{displaymath}
\inf_{\Omega} \big( u( t, u_{0},f_{0}) -u(t, u_{1}, f_{1})\big)>0  ,\; \mbox{ for all }\;t\ge0,  
\end{displaymath}
we say that (\ref{nonlinear_diffusion}) satisfies a strong  
comparison principle.
  
\end{definition}

With respect to maximum principles we define the following. 

\begin{definition} 
For $u_{0} \in X$ denote by $u(t,u_{0}, f)$ the solution of
(\ref{nonlinear_diffusion}) which we assume to exist for some class
of nonlinear terms $f$.

\noindent i) 
We say that (\ref{nonlinear_diffusion}) satisfies a weak
maximum principle if for any  
$u_0 \in X$ such that 	$u_0\ge 0$, then
\begin{displaymath}
u( t, u_{0},f)\ge 0,\; \mbox{ for all }\;t\ge0. 
\end{displaymath}

\noindent ii) If moreover $u_{0} \neq 0$ 
and 
\begin{displaymath}
u( t, u_{0},f)> 0,\; \mbox{ for all }\;t\ge0,  
\end{displaymath}
we say that (\ref{nonlinear_diffusion}) satisfies a strict 
maximum principle.

\noindent iii) If furthermore 
\begin{displaymath}
\inf_{\Omega} u( t, u_{0},f) >0  ,\; \mbox{ for all }\;t\ge0,  
\end{displaymath}
we say that (\ref{nonlinear_diffusion}) satisfies a strong  
maximum  principle.
  
\end{definition}

\subsection{The case of  a globally Lipschitz reaction term}
\label{exist_uniq_pos_Glob_Lipschitz}

In this section we will assume $f: \Omega \times \R \to \R$ is
globally Lipschitz in the second variable. This implies that for either
$X=L^{p}(\Omega)$, $1\leq p\leq \infty$ or $X=C_{b}(\Omega)$, the
Nemitcky operator  $F:X\rightarrow X$ is
Lipschitz. 

%

We can even consider here a little more general situation by
considering a  general globally Lipschitz operator $G:X\to X$, not
necessarily a Nemitcky operator,  and 
consider  the nonlocal nonlinear problem 
\begin{equation}\label{eq:nonlinear_global_lipchitz}
\left\{\!\!\!
\begin{array}{lll}
u_t(x,t)\!\!\!\!\! & =Lu(x,t)+G(u)(x,t),\!\! & x\in\Omega, t\in \R,\\
u(x,0)\!\!\!\!\! & =u_0(x),\!\! & x\in\Omega,
\end{array}
\right.
\end{equation}
where, as above, $L=K-hI$. Observe that for a Nemitcky operator
(\ref{eq:Nemitcky}), if $X =L^{p}(\Omega)$, $1\leq p\leq \infty$ then
we will require $f(x,s)$ to be Lipchitz in $s\in \R$, uniformly in
$x\in \Omega$ and $g(x)= f(x,0) \in X$. On the other hand, if
$X=C_{b}(\Omega)$ then we will additionally requiere $f(x,s)$
continuous in $(x,s)\in \Omega \times \R$.

The following result gives the existence and uniqueness of solutions
to (\ref{eq:nonlinear_global_lipchitz}).

\begin{proposition}\label{existence_global_lipschitz_G}

Assume, as in
Proposition \ref{positive_solution_with_u_0_positive}, 
$J\in L^{p}(\Omega, L^{p'}(\Omega))$ for some $1\leq p \leq \infty$ and then denote $X=
L^{p}(\Omega)$ and assume $h\in L^{\infty}(\Omega)$, or $J\in
C_{b}(\Omega, L^{1}(\Omega))$ and then denote $X= C_b(\Omega)$ and  assume $h\in C_b(\Omega)$.

Then problem  (\ref{eq:nonlinear_global_lipchitz}) 
	has a unique global solution $u\in \mathcal{C}((-\infty,\infty), X)$, for 
	every $u_0\in X$, given by 
	\begin{equation}\label{nonlinear_VCF_G_Lipschitz}
		u(\cdot,t)=e^{Lt}u_0+\displaystyle\int_0^t e^{L(t-s)}G(u)(\cdot,s)
		\,ds, \qquad t\in \R.
	\end{equation}
	Moreover, $u\in \mathcal{C}^1((-\infty,\infty), X)$  is a strong 
	solution of (\ref{eq:nonlinear_global_lipchitz}) in 	$X$.
\end{proposition}
 This result can be proved using a fixed point argument using the
 variation of constants formula in $\mathcal{C}([-\tau, \tau], X)$
 for some $\tau>0$ independent of the initial data, and a prolongation
 argument. As the arguments are standard,  we will
 omit the proof. Also note that the variation of constans formula,
 that is, the right hand side of (\ref{nonlinear_VCF_G_Lipschitz}) maps
 $L^{1}([-\tau, \tau], X)$ into $\mathcal{C}([-\tau, \tau], X)$. Hence
 the solution of (\ref{eq:nonlinear_global_lipchitz}) is unique in both
 spaces. The fact that (\ref{nonlinear_VCF_G_Lipschitz}) is a strong solution of
 (\ref{eq:nonlinear_global_lipchitz}) follows from  Theorem in
 \cite[p. 109]{Pazy}.

 \begin{remark}\label{rem:alt_formulation_VCF}
   Observe that for any $\beta \in \R$ we can  rewrite  (\ref{eq:nonlinear_global_lipchitz}) as  
   \begin{equation} \label{eq:alt_nonlinear_global_lipchitz }
u_t(x,t)  =Lu(x,t)-\beta\, u(x,t)+G(u)(x,t)+\beta\, u(x,t) . 
   \end{equation}
Since 	$L-\beta I$ and $G + \beta I$ satisfy the same assumptions as
in Proposition \ref{existence_global_lipschitz_G} then we obtain the
alternative reperesentation of the solution of
(\ref{eq:nonlinear_global_lipchitz}) as 
\begin{equation}
  \label{eq:alt_nonlinear_VCF_G-Lipschitz}
u(t)= e^{(L-\beta I)t}u_{0}+ \int_0^t e^
		{(L-\beta I)(t-s)}\left(G(u)(s)+\beta u(s)\right) ds,  \qquad t\in \R
\end{equation}
and Proposition \ref{existence_global_lipschitz_G} remains true under
this alternative formulation. 

 \end{remark}

Now we  prove monotonicity properties with respect to the 
nonlinear term  for  problem
(\ref{eq:nonlinear_global_lipchitz}). Notice that for  the case of the
Nemitcky operator (\ref{eq:Nemitcky}) if $f$ is globally Lipschitz in
the second variable, then  there exists a constant
  $\beta>0$, such that $u\mapsto f(x,u)+\beta u$ is increasing, for
  all $x\in \Omega$. Hence the   assumptions below on the nonlinear terms are 
satisfied.


\begin{proposition}\label{comparison_G1G2_Lipschitz}
{\bf (Weak, strict  and strong comparison principles) }
  

  Under the assumptions in Proposition
  \ref{existence_global_lipschitz_G},  consider globally Lipschitz
  functions  $G:X\rightarrow X$ such that there exists a constant
  $\beta>0$, such that $G+\beta I$ is increasing. Then 

\noindent i) If $J\geq 0$ then  problem (\ref{nonlinear_VCF_G_Lipschitz}) satisfies a weak
  comparison principle. 
	
\noindent ii) If, additionally,  $\Omega$ is $R$-connected and $J$ satisfies 
hypothesis (\ref{eq:J_positive}) then problem
(\ref{nonlinear_VCF_G_Lipschitz}) satisfies a strict 
  comparison principle. 


\noindent iii) Finally, if  moreover  $|\Omega|<\infty$, $\Omega$ is  $R$--connected and  the measure satisfies
  (\ref{eq:measure_non_degenerates}) and $J$ satisfies
  (\ref{eq:J_strict_positive})  then problem
(\ref{nonlinear_VCF_G_Lipschitz}) satisfies a strong 
  comparison principle. 


\end{proposition}
\begin{proof}
Take  two such functions such that $G_{0} \geq G_{1}$ and  take
$\beta$ such that both $G_{i}+\beta I$, $i=0,1$, are
increasing. Denote 
$u^i(t) = u( t, u_{i},G_{i})$, $t\in \R$ the corresponding solutions
of (\ref{nonlinear_VCF_G_Lipschitz}), which by Remark
\ref{rem:alt_formulation_VCF} and Proposition
\ref{existence_global_lipschitz_G} are the unique fixed points of 
	\begin{displaymath} 
		\mathcal{F}_i(u)(t)=e^{(L-\beta I)t}u_i+\displaystyle\int_0^t e^{(L-
		\beta I)(t-s)}\left(G_i(u)(s)+\beta u(s)\right) ds
	\end{displaymath}
	in $V=\mathcal{C}([-\tau,\tau],X)$ because $\mathcal{F}_i$ is
        a contraction in $V$  provided $\tau$ small enough, for $i=0,\, 
	1$ (independent of the initial data). Consider then the
        sequence of Picard iterations,
        $u^i_{n+1}(t)=\mathcal{F}_i(u^i_n)(t)$,   $n=1, 2,
           \ldots$,  $0\leq t \leq \tau$,  $u^{i}_{1} = u_{i}$.
	 Then the sequence $u^i_n(\cdot,t)$ converges to $u^i(\cdot,t)$ in $V$. 
	 Now, we 
	 are going to prove that the solutions are ordered for all $t\ge 0$. 
	 We take the first term of the Picard iteration as
         $u^0_1(x,t)=u_0(x) \geq u_1(x) = u^1_1(x,t)$, 
	 then 
	 \[
	 u^i_2(t)=\mathcal{F}_i(u^i_1)(t)=e^{(L-\beta I)t}u_i +\displaystyle
	 \int_0^t e^
	 {(L-\beta I)(t-s)}\left(G_i(u_i)+\beta u_i\right) ds ,  \quad 0\leq t \leq \tau .
	 \]

\noindent i)  If $J\geq 0$, by Proposition
\ref{positive_solution_with_u_0_positive}, $e^{(L-\beta I)t}u_0\ge
e^{(L-\beta I)t}u_1$  for $t\in [0,\tau]$ 
and since
         $G_0+\beta I\ge G_1+\beta I$ and are increasing, we have 
$$
e^{(L-\beta I)(t-s)}\!\left(G_0(u_0)\!+\!\beta u_0\right)\!\ge\! e^{(L-\beta I)(t-
	 s)}\!\left(G_1	 (u_1)\!+\!\beta u_1\right),\, \quad  0 \leq s
       \leq t \leq \tau.
$$
	 Hence $u^0_2(t)\ge u^1_2(t)$ for all $t \in[0,\tau]$.  
	 Repeating this argument, we obtain that $u^1_n(t)\ge
         u^2_n(t)$   for all  $t\in [0,\tau]$  for    every $n\ge1$. 
	  Since $u^i_n(t)$ converges to $u^i(t)$, in $V$, then
          $u^0(t)\ge u^1(t)$  for $ t\in [0,\tau]$.
	 
	  Now, we consider the solutions of
          (\ref{nonlinear_VCF_G_Lipschitz}) 
	  with initial data $u^0(\tau) \geq u^1(\tau)$, and   arguing
          as above  we obtain that 
	  $u^0(t)\ge u^1(t)$ for all $t\in [\tau,2\tau]$.  Repeating this argument, we obtain that 
	  \begin{equation}\label{ineq_g_1_g_2_eq}
	  	  u^0( t)\ge u^1( t),\; \mbox{ for all }\;t\ge0.
	  \end{equation}
	  
\noindent ii) Using (\ref{ineq_g_1_g_2_eq}),
$u^i(t)=\mathcal{F}_i(u^i)(t)$ 
and that   $G_0+\beta I\ge G_1+\beta I$ and are increasing, we get 
	  \begin{displaymath}\label{ineq_g_1_g_2_eq_2}
	  G_0(u^0)(t)+\beta u^0(t) \ge G_0(u^1)(t)+\beta u^1(s)
	  \geq G_1(u^1)(t)+\beta u^1(t), \quad  t\ge0.
	  \end{displaymath}
	  From 	  Proposition
          \ref{positive_solution_with_u_0_positive}, we have that
          $e^{(L-\beta I)t}\left(G_0(u^0)(s)+\beta u^0(s)\right) >
          e^{(L-\beta I)t}\left(G_1 
	  (u^1)(s)+\beta u^1(s)\right)$ for   $0\leq s \leq t$. 
	  Therefore, 
	   \begin{displaymath}
	   \int_0^t \!e^{(L-\beta I)(t-s)}\big( G_0(u^0)(s)+\beta
           u^0(s)\big) > \int_0^t e^{(L-\beta I)(t-s)} \left(G_1(u^1)(s)+\beta u^1(s)\right) ds,
	  \end{displaymath}
	  for all $t> 0$. Thus, $u^1(t)>u^2(t)$, for all $t>0$.

\noindent iii)
In this case, by Proposition \ref{prop:Linear_Parab_Strong_MaxPple},
$\inf_{\Omega} e^{(L-\beta  I)t} (u_0 -u_{1}) >0$,  $t>0$ 
and we get the result.           
\end{proof}



Concerning maximum principles, we get the following result. 

\begin{proposition}\label{nonlinear_max_pple_G_Lipschitz}
{\bf (Weak, strict  and strong parabolic maximum principle) }

Under the assumptions in Proposition \ref{existence_global_lipschitz_G}, assume 
	$G:X\rightarrow X$ globally Lipschitz,  there exists a constant 
	$\beta>0$, such that $G+\beta I$  is increasing and 	$G(0)\ge 0$. 

\noindent i) If $J\geq 0$ then  problem (\ref{nonlinear_VCF_G_Lipschitz}) satisfies a weak
  maximum  principle. 
	
\noindent ii) If, additionally,  $\Omega$ is $R$-connected and $J$ satisfies 
hypothesis (\ref{eq:J_positive}) then problem
(\ref{nonlinear_VCF_G_Lipschitz}) satisfies a strict 
  maximum  principle. 

\noindent iii) Finally, if  moreover  $|\Omega|<\infty$, $\Omega$ is  $R$--connected and  the measure satisfies
  (\ref{eq:measure_non_degenerates}) and $J$ satisfies
  (\ref{eq:J_strict_positive})  then problem
(\ref{nonlinear_VCF_G_Lipschitz}) satisfies a strong 
maximum  principle.


\end{proposition}
\begin{proof}
We have that $u(t) = u(t,u_{0}, G)$, $t\geq 0$, is 
	the unique fixed point of 
	\begin{displaymath} 
	\mathcal{F}(u)(t)=e^{(L-\beta I)t}u_0+\displaystyle\int_0^t e^{(L-\beta I)
	(t-s)}\left(G(u)(\cdot,s)+\beta u(s)\right) ds
	\end{displaymath}
	which  is a contraction in $V=\mathcal{C}([-\tau,\tau],X)$, for
        some $\tau$  small  
	enough independent of the initial data.		We consider
        the sequence of Picard iterations,
        $u_{n+1}(t)=\mathcal{F}(u_n)(t)$,   $n=1,2 \ldots$, $0\le t\le \tau$
which   converges to $u(\cdot,t)$ in         $V$.

         \noindent i)
         We take  
	 $u_1(t)=u_0\geq 0$, then 
	 \[
	 u_2(t)=\mathcal{F}(u_1)(t)=e^{(L-\beta I)t}u_0+\displaystyle\int_0^t e^{(L-\beta I)
	 (t-s)}\left(G(u_0)+\beta u_0\right) ds.
	 \] 
	 Thanks to Proposition
         \ref{positive_solution_with_u_0_positive} we have
         $e^{(L-\beta I)t}u_0\ge 0$ for  $ t\in [0,\tau]$ 
while on  the other hand, since  $G(0)\ge 0$,  $\beta>0$ 
	 and $\;G(\cdot)+ \beta I$ is increasing,  then $G(u)+\beta u\ge 0$ for 
	 all $u\ge 0$. Hence, from Proposition 
	 \ref{positive_solution_with_u_0_positive} we obtain that
         $e^{(L-\beta I)(t-s)}\left(G(u_0)+\beta u_0\right)\ge 0$, for all $0 \leq s \leq t \leq \tau$. 
	 Hence, 
         $u_2(t)\ge 0$ for all  $t\in[0,\tau]$. 

	 Repeating this argument, we get that $u_n(t) \geq 0$ for
         every $n\ge1$ and $t \in [0,\tau]$. Since $u_n(t)$ converges
         to $u(t)$ in $V$ then  $u(t)\geq 0$ 	 for all $t$ in $[0,\tau]$.
	 
Now we consider the solution of (\ref{nonlinear_VCF_G_Lipschitz}) with
initial data $u(\tau)\geq 0$ and  arguing as above we have that
$u(t)\geq 0$ is nonnegative for all $t\in[\tau,\,2\tau]$ and  thus for
$t\in [0,2\tau]$. 
           Repeating this argument, we prove that $u(t) \geq 0$  for all $t\ge 0$.

           \noindent ii)   Using part i), $u(t)=\mathcal{F}(u)(t)$, 
Proposition \ref{positive_solution_with_u_0_positive}  and 
	 $\big(G+\beta I\big)(u)\ge 0$ for all $u\ge 0$, we get that
         $e^{(L-\beta I)t}u_0>0$,  for $t>0$
and $\int_0^t e^{(L-\beta I)(t-s)}\left(G(u)(x,s)+\beta u(x,s)\right)
ds\ge 0$,  for $t\ge 0$.
	 Thus, 
         we 	 have that $u(t)>0$ for all $t>0$.

\noindent iii) In this case, from Proposition
\ref{prop:Linear_Parab_Strong_MaxPple} we have $\inf_{\Omega}
e^{(L-\beta I)t} u_0  >0$  for $t>0$, 
 and we conclude. 
\end{proof}
%


Now we  introduce the definition of supersolution and subsolution  to 
(\ref{eq:nonlinear_global_lipchitz}).

\begin{definition} 
	Let $X=L^p(\Omega)$, with $1\le p\le\infty$ or  $X=\mathcal{C}_b
	(\Omega)$, we say that $\overline{u}\in \mathcal{C}([a,b], X)$ is a {\bf 
	supersolution} 
	to  (\ref{eq:nonlinear_global_lipchitz}) in $[a,b]$, if for any $t\ge s$, with 
	$s,t\in[a,\,b]$
	\begin{equation}\label{def_supersol}
		\displaystyle \overline{u}(t)\ge e^{L(t-s)}\overline{u}(s)+
		\int_s^t e^{L(t-r)}G(\overline{u})(r)dr.
	\end{equation}
	We say that $\underline{u}$ is a {\bf subsolution} if the reverse 
	inequality holds.
\end{definition}

\begin{remark} \label{rem:supersolutions}

\noindent i) As above, using (\ref{eq:alt_nonlinear_global_lipchitz
}) and (\ref{eq:alt_nonlinear_VCF_G-Lipschitz}) we have the following
alternative defnition of supersolutions of  (\ref{eq:nonlinear_global_lipchitz}) in $[a,b]$, if for any $t\ge s$, with 
	$s,t\in[a,\,b]$
	\begin{equation}\label{alt_def_supersol}
		\displaystyle \overline{u}(t)\ge e^{(L-\beta I) (t-s)} \overline{u}(s)
                + \int_s^t e^
		{(L-\beta I)(t-r)} \left(G(\overline{u})(r)+\beta \overline{u}(r)\right) dr
	\end{equation}
with an analogous definition of  {\bf subsolution} with  reverse
inequality.

\noindent ii) 
Assuming  $J\geq 0$ if $\overline{u}\in \mathcal{C}([a,b], X)$
        is 	differentiable and satisfies that 
	\begin{displaymath} 
		\overline{u}_t(t)\ge L \overline{u} (t)+G(\overline{u})(t),\;\mbox
		{ for }\; t\in[a,\,b]
	\end{displaymath}
	then $\overline{u}$ is a supersolution in the sense of
        (\ref{def_supersol}) or (\ref{alt_def_supersol}). The same happens for subsolutions if the
        reverse  	inequality holds. 
%
%
%

\end{remark}

The following proposition states  that supersolutions and subsolutions of
(\ref{eq:nonlinear_global_lipchitz})  are above and below solutions respectively.

\begin{proposition}
  \label{ordered_solutions_supersolutions_G_Lipschitz}

Under the assumptions in Proposition
\ref{existence_global_lipschitz_G}, assume $J\geq 0$ and 
	$G:X\rightarrow X$ is globally Lipschitz and there exists a constant 
	$\beta>0$, such that $G+\beta I$  is increasing. 
	Let $u( t,\, u_0)$ be the solution to
        (\ref{eq:nonlinear_global_lipchitz}) with initial data
        $u_0\in X$, and let $\bar{u}(t) $ be a supersolution to
        (\ref{eq:nonlinear_global_lipchitz}) in $[0,T]$.
	
	If  $\overline{u}(0)\ge u_0$, then 
	$$\bar{u}(t)\ge u( t,\, u_0),\quad\mbox{ for }\; t\in [0,T].$$
	The same is true for subsolutions with reversed inequality.
\end{proposition}
\begin{proof}
We have that $u(t) = u(t,u_{0},G)$, $t\geq 0$,  	the unique fixed point of 
	\begin{displaymath} 
	\mathcal{F}(u)(t)=e^{(L-\beta I)t}u_0+\displaystyle\int_0^t e^{(L-\beta I)
	(t-s)}\left(G(u)(\cdot,s)+\beta u(s)\right) \,ds \quad t\in
      [0,\tau] 
	\end{displaymath}
	in $\mathcal{C}([0,\tau],X)$,  for some $\tau$  small  enough
        independent of the initial data. Also, we can assume without
        loss of generality that $\tau \leq T$.  
Also, we  consider the sequence of Picard iterations in $V=\mathcal{C}
	([0,\tau],X)$, $u_{n+1}(t)=\mathcal{F}(u_n)(t)$,
        $n=1,2\ldots$, $ t\in  [0,\tau] $ 
	 with $u_1(t)=\overline{u}(t)$. Then the sequence $u_n(t)$ 
	 converges to $u(t)$ in $V$ and we show below that $\bar{u}
         (t) \ge u_{n}(t)$,  for  $n=1,2,\ldots$ and $ t\in
         [0,\tau]$.  


Note that $u_1= \bar{u}$ satisfies by definition $\bar{u}(t)\ge
\mathcal{F}(\bar{u})(t)$ for $t\in[0,\tau]$ and then 
         we have that $\bar{u}(t)\ge \mathcal{F}(\bar{u})(t)=u_2(t)$,  $t\in [0, \tau]$. 
Observe now that from the proof of Proposition
\ref{comparison_G1G2_Lipschitz} we have that $\mathcal{F}$ is
increasing in $V$, and therefore $\overline{u} (t)\ge \mathcal{F}(\overline{u})(t) \ge \mathcal{F}
		(u_2)(t)=u_{3}(t)$,   $t\in [0, \tau]$. 
By induction we get the claim.  Since 
	$u_n(t)$ converges to $u(t)$ in $V$ we have that   $\bar{u}(t)\ge u(t,u_0)$, $t\in [0, \tau]$. 
Repeating this argument with initial data $u(\tau)\leq \bar u(\tau)$
we get the result in $[0,T]$.  
      \end{proof}

\subsection{The case of a  locally Lipchitz reaction term} 
\label{existence_sol_with_f_loc_Lip}


In this section we consider (\ref{nonlinear_diffusion}) with some
classes of nonlinear functions  $f: \Omega \times \R \to \R$. More
precisely we consider 
\begin{equation}\label{nonlinear_loc_lip}
	\!\!\!\left\{\!\!\!\!
	\begin{array}{rcl}
	u_t(x,t) \!\!\!\!& =&\!\!\!\! Lu(x,t)\!+\!f(x,u(x,t)),\quad x\in\Omega,\, t>0,\\
	u(x,t_0) \!\!\!\! & =& \!\!\!\! u_0(x),\quad x\in\Omega, 
	\end{array}
	\right.
\end{equation}
 with $L=K-hI$ and $f:\Omega\times\mathbb{R}\rightarrow \mathbb{R} $ 
such that $f(x,s)$ is locally Lipschitz in the 
 variable $s\in\mathbb{R}$, uniformly with 
respect to $x\in\Omega$, and $f$ satisfies the sign condition 
(\ref{hyp_f_Cs_2_D}) below. If we work in 
$X=C_{b}(\Omega)$ then we will additionally requiere $f(x,s)$
continuous in $(x,s)\in \Omega \times \R$.

The strategy we use to solve (\ref{nonlinear_loc_lip}) is as follows.
For $k>0$, let us introduce a globally lipschitz function,
$f_k:\Omega\times\R\to\R$,  such that
		\begin{equation}\label{efe_k} 
			f_k(x,u)=f(x,u)\qquad \mbox{for }\; |u|\le k, \;  \mbox{ and  }\; 
			x\in \Omega.
		\end{equation}
Then we consider the problem 
\begin{equation}\label{nonlinear_efe_k}
\left\{\!\!\!
\begin{array}{lll}
u_t(x,t)\!\!\!\!\! & = Lu(x,t)+F_k(u)(x,t),\quad x\in\Omega,\, t>0,
\\
u(x,0)\!\!\!\!\! & =u_0(x), \quad x\in\Omega,\!\! & 
\end{array}
\right.
\end{equation}
where $F_k:X\to X$ is the globally Lipschitz  Netmitcky operator associated to the globally 
Lipschitz function $f_k$. Then Proposition
\ref{existence_global_lipschitz_G} gives the  existence and uniqueness of solutions to  
(\ref{nonlinear_efe_k}). 


Also, for fixed $k$, since $f_k$ is globally Lipschitz,  there exists $\beta>0$ such 
that $f_k+\beta I$ is increasing, 
and then $F_k+\beta I$ is increasing. Hence, we can apply Propositions 
\ref{comparison_G1G2_Lipschitz}, 
\ref{nonlinear_max_pple_G_Lipschitz} and 
\ref{ordered_solutions_supersolutions_G_Lipschitz} 
for the problem (\ref{nonlinear_efe_k}).

Hence, in order to solve (\ref{nonlinear_loc_lip}) we 
show that under the sign condition (\ref{hyp_f_Cs_2_D}) on $f$,  for bounded initial data 
we estimate the sup norm of the solution
$u_{k}(t)$, see Proposition
\ref{global_sol_initial_data_bounded_2}. Later we solve
(\ref{nonlinear_loc_lip}) for initial data in $X=L^{p}(\Omega)$ by
assuming some natural growth condition on $f$, see Theorem
\ref{global_solutions_f_growth} below.

\begin{proposition}\label{global_sol_initial_data_bounded_2}

Assume $0\leq J\in L^{\infty}(\Omega, L^{1}(\Omega))$ and then denote $X=
L^{\infty}(\Omega)$ and assume $h\in L^{\infty}(\Omega)$, or $J\in
C_{b}(\Omega, L^{1}(\Omega))$ and then denote $X= C_b(\Omega)$ and
assume $h\in C_b(\Omega)$. Also, assume that  the locally Lipschitz 
	function  $f: \Omega \times \R \to \R$  satisfies that
        $g=f(\cdot,0) \in X$ and there exist $C,\,D\in\mathbb{R}
	$, with $C>0$ and $D\ge 0$ such that 
	\begin{equation}\label{hyp_f_Cs_2_D}
	f(\cdot,s)s\le Cs^2+D|s|,\;\quad s\in \R.
	\end{equation}
	
Then there exists a unique global solution of (\ref{nonlinear_loc_lip}) with 
	initial data $u_0 \in X$, such that $u(\cdot,t)$ in $\mathcal{C}([0,T], 
	X)$, for all $T>0$, with
	\begin{equation}\label{var_const_form_exist_mejorada}
		u(\cdot,t)=e^{Lt}u_0+\displaystyle\int_0^t e^{L(t-s)}f(\cdot,u
		(\cdot,s))\,ds, \quad t\geq 0.
	\end{equation}
	Moreover, we have that  $u$ is a strong solution of (\ref{nonlinear_loc_lip}) in 	$X$.
\end{proposition}
\begin{proof}
First of all, observe that the function  $h_{0}$ defined in
(\ref{eq:definition_h0}),  belongs to $X$ and 
let us prove that $(h_0-h)s+f(\cdot,s)$ satisfies   
	hypothesis (\ref{hyp_f_Cs_2_D}). 	
	Since $f$ satisfies (\ref{hyp_f_Cs_2_D}) and $h,h_0\in X$, then
	\begin{equation}\label{hyp_para_s_con_h_y_h_0}
	\begin{array}{lll}
	(h_0-h)s^2+f(\cdot,s)s 
	& \le \big(\|h_0-h\|_{L^{\infty}(\Omega)}
	+C\big)s^2+D|s|
	& \le C_1s^2+D|s|.
	\end{array}
	\end{equation}
We will denote $C_{1}$ again by $C$ in order to  to simplify the notation. 

	Fix $0<M\in\mathbb{R}$. We introduce the auxiliary problem $	\left\{
	\begin{array}{ll}
	\dot{z}(t) & =Cz(t)+D\\
	z(0) & =M.
	\end{array}
	\right.$ 
Since $z(t)$ is increasing, for any given  $T>0$ 
	\begin{equation}\label{cota_z_mejorada}
	0\le z(t)\le z(T)\quad  t\in
	[0,T].
	\end{equation}	
	Given $T>0$ and $M>0$, from (\ref{cota_z_mejorada}) we choose
        $k\ge z(T)$ 
and consider a globally Lipchitz truncation of $f$,  $f_{k}$, as in
(\ref{efe_k}). Thus 
	\begin{equation}\label{f_igual_f_k_mejorada}
		f_k(\cdot,z(t))=f(\cdot,z(t)),\quad  t\in [0,T].
	\end{equation}

	We prove below that $z$ is a supersolution of
        (\ref{nonlinear_efe_k}) in  $[0,T]$. First since $z(t)$ is independent of the 
	variable $x$, we have that  $Kz(t)=h_0z(t)$. Now since  $z(t)$ is nonnegative
        for all $t\in[0,T]$, using  	(\ref{f_igual_f_k_mejorada})
        and (\ref{hyp_para_s_con_h_y_h_0}) we get,  for $t\in[0,T]$, 
	\[
	Kz(t)-hz(t)+f_k(\cdot,z(t))  = (h_0-h) z(t)+f_k
	(\cdot,z(t))
	 \le  Cz(t)+D=\dot{z}(t) . 
	\]
	Hence, $z$ is a supersolution of (\ref{nonlinear_efe_k}) 
	in  $[0,T]$, see Remark \ref{rem:supersolutions}.
	
	Analogously, let us consider the auxiliary problem $	\left\{
	\begin{array}{ll}
	\dot{w}(t) & =Cw(t)-D\\
	w(0) & =-M.
	\end{array}
	\right.$ 
	Then $w(t)=-z(t)$, and we obtain that  
%
        $|w(t)|<z(T)$ for $ t\in[0,T]$. 
	Moreover $w$ is a subsolution of (\ref{nonlinear_efe_k}) in  $t\in
	[0,T]$. 
Hence, if  $\|u_0\|_X\le M$, from Proposition
\ref{ordered_solutions_supersolutions_G_Lipschitz}, we obtain
$w(t)\le u_k(t,u_0)\le z(t)$, for  $t\in [0,T] $  
and therefore 
 $$
|u_k(t,u_0)|\le z(T)\le k\;\mbox{ for all }\;t\in [0,T].
$$
In particular $f_k(\cdot, u_k(t,u_0))=f(\cdot, u_k (t,u_0))$ and thus,
$u_k(x,t,u_0)$ is a (strong) solution to (\ref{nonlinear_loc_lip}) in
$t\in [0,T]$ and satisfies (\ref{var_const_form_exist_mejorada}).

	
	Now, let us prove uniqueness.   
	Consider a solution $u\in
	 \mathcal{C}
	([0,T],X)$  of  problem (\ref{nonlinear_loc_lip}) with initial data 
	$u_0\in X$, given by (\ref{var_const_form_exist_mejorada}). Since  
	$u\in \mathcal{C}([0,T],X)$, then 
	$$\sup\limits_{t\in [0,T]}\sup_{x\in\Omega}|u(x,t,u_0)|<\widetilde{C}.$$ 
	Thus, if we choose $k>\widetilde{C}$, then 
	$f_k(\cdot, u(\cdot, t))=f(\cdot, u(\cdot, t))$
	and then $u$ coincides  in $[0,T]$ with the
        solution of (\ref{nonlinear_efe_k}), $u_k$. 
        Thus, we have  the uniqueness 
	of the solution of (\ref{nonlinear_loc_lip}).
\end{proof}

\begin{remark} 
Observe that 	hypothesis (\ref{hyp_f_Cs_2_D}) on $f$ in   Proposition 
	 \ref{global_sol_initial_data_bounded_2}  
	is somehow  optimal. 
For this  assume $J(x,y)=J(y,x)$, and  consider the problem
	\begin{equation}\label{eq_Kaplan}
	\left\{
	\begin{array}{l}
	u_t (x,t)  \displaystyle =\int_{\Omega}J(x, y)u(y,t)\, dy-h(x)u(x,t)+u^{\rho}(x,t)\\
	u(x,0)=u_0(x) 
	\end{array}
	\right.
	\end{equation}
	with $\rho >1$,  $u_0\in L^{\infty}(\Omega)$, and $u_0\ge
        0$.

        Observe first that the operator $K$ has a principal eigenvalue
        $\Lambda$, see  (\ref{eq:oscilation_h}) with $h=0$ cf.
        \cite{RB2017_max_pples}. 
	Let $\phi>0$ be an eigenfunction associated to  $\Lambda$ normalized as $\int_{\Omega}\phi(x)dx=1$ and define
        $z(t)=\int_{\Omega}u(t)\phi$. Then 
	\begin{equation*} 
	\begin{array}{rcl}
		\displaystyle{dz \over dt}(t)  &=& \displaystyle\int_{\Omega}u_t(x,t)\phi(x)dx
		\smallskip\\
		&=& \displaystyle \int_{\Omega}\!\int_{\Omega}\!\!J(x, y)\phi(x)dx\, 
		u(y,t)dy
		-\!\!\int_{\Omega}\!\!h(x)\phi(x)u(x,t)dx+\!\!\int_{\Omega}\!\!u^{\rho}(x)\phi(x)dx.
	\end{array}
	\end{equation*}
	Since $J(x,y)=J(y,x)$ and $\phi$ is an eigenfunction of $K$ we have that 
	\begin{equation*} 
	\begin{array}{ll}
		\displaystyle{dz \over dt}(t) &= \displaystyle \int_{\Omega}\int_{\Omega}
		\!J(y,x)\phi(x)\, dx \, u(y,t)dy
		-\!\!\int_{\Omega}\!h(x)\phi(x)u(x,t)dx+\!\int_{\Omega}\!u^{\rho}
                                                                  (x,t)\phi(x)dx
		\smallskip\\
		&= \displaystyle \Lambda \int_{\Omega} \phi(y)u(y,t)\,
                  dy  	-\!\!\int_{\Omega}\!h(x)\phi(x)u(x,t)dx + \int_{\Omega}u^{\rho}(x,t)\phi(x)dx
		\smallskip\\
		&= \displaystyle \Lambda z(t)
                  -\!\!\int_{\Omega}\!h(x)\phi(x)u(x,t)dx
                  +\int_{\Omega}u^{\rho} (x,t)\phi(x)dx.
	\end{array}
	\end{equation*}
Then using Jenses's inequality we get 
	\begin{equation*} 
{dz \over dt}(t) \geq \big( \Lambda -\|h\|_{\infty} \big) z(t)  +
z^{\rho} (t) = F(z(t)).
	\end{equation*}
        Thus, if $z(0)$ is sufficiently large the solution of
        (\ref{eq_Kaplan}) must cease to exist in finite time. 
\end{remark}

Since we have proved that the solutions of (\ref{nonlinear_loc_lip})
for initial data  $u_0$ in $L^{\infty}(\Omega)$ or in
$\mathcal{C}_b(\Omega)$ coincides,  on a given time interval,  with a solution of some  problem
(\ref{nonlinear_efe_k}), with a truncated globally Lipschitz function
$f_k$, these solutions inherit all the monotonicity properties in
Section   \ref{exist_uniq_pos_Glob_Lipschitz} as we now state.



 
\begin{corollary}\label{comparison_f1f2_Loc_Lipschitz}
  {\bf  (Weak, strict  and strong comparison principles)}
  
	Under the hypotheses of  Proposition
        \ref{global_sol_initial_data_bounded_2},
	then for any initial data $u_0\in X =L^{\infty}(\Omega)$ or
        $X=\mathcal{C}_b(\Omega)$

\noindent i) If $J\geq 0$ then  problem (\ref{nonlinear_loc_lip}) satisfies a weak
  comparison principle. 
	
\noindent ii) If, additionally,  $\Omega$ is $R$-connected and $J$ satisfies 
hypothesis (\ref{eq:J_positive}) then problem (\ref{nonlinear_loc_lip}) satisfies a strict 
  comparison principle. 

\noindent iii) Finally, if  moreover  $|\Omega|<\infty$, $\Omega$ is  $R$--connected and  the measure satisfies
  (\ref{eq:measure_non_degenerates}) and $J$ satisfies
  (\ref{eq:J_strict_positive})  then problem (\ref{nonlinear_loc_lip}) satisfies a strong 
  comparison principle.

	

\end{corollary}

\begin{corollary} \label{nonlinear_max_pple_f_Loc_Lipschitz}
{\bf (Weak and strong parabolic maximum principle)}

Under the hypotheses of  Proposition
        \ref{global_sol_initial_data_bounded_2}, 
        assume moreover that  $f(\cdot, 0)\ge 0$. 
        Then for initial data in $X =L^{\infty}(\Omega)$ or
        $X=\mathcal{C}_b(\Omega)$,

\noindent i) If $J\geq 0$ then  problem (\ref{nonlinear_VCF_G_Lipschitz}) satisfies a weak
  maximum  principle. 
	
\noindent ii) If, additionally,  $\Omega$ is $R$-connected and $J$ satisfies 
hypothesis (\ref{eq:J_positive}) then problem
(\ref{nonlinear_VCF_G_Lipschitz}) satisfies a strict 
  maximum  principle. 

\noindent iii) Finally, if  moreover  $|\Omega|<\infty$, $\Omega$ is  $R$--connected and  the measure satisfies
  (\ref{eq:measure_non_degenerates}) and $J$ satisfies
  (\ref{eq:J_strict_positive})  then problem
(\ref{nonlinear_VCF_G_Lipschitz}) satisfies a strong 
maximum  principle.


\end{corollary}

\begin{corollary} \label{ordered_solutions_supersolutions_f_Loc_Lipschitz}

Under the hypotheses of  Proposition
        \ref{global_sol_initial_data_bounded_2}, 
        let   $u( t,\, u_0)$ be a solution to (\ref{nonlinear_loc_lip}) with initial 
	data $u_0\in X =L^{\infty}(\Omega)$ or
        $X=\mathcal{C}_b(\Omega)$, and let $\bar{u}( t)$ be a supersolution to 
	(\ref{nonlinear_loc_lip}) in $[0, T]$. 
	
	If $\bar{u}(0)\ge u_0$, then 
	$$\bar{u}( t)\ge u( t,\, u_0), \quad \mbox{ for all }\; t\in [0, T].$$
	The same is true for subsolutions with reversed inequality.
\end{corollary}

Now we prove the existence and uniqueness for the problem (\ref{nonlinear_loc_lip}) with initial
data in $L^q(\Omega)$ for some suitable $1< q< \infty$, depending on
the growth of $f$. 


\begin{theorem}\label{global_solutions_f_growth}
	Assume  $|\Omega|<\infty$, $J(x,y)=J(y,x)$ and 
$0\leq J\in L^{\infty}(\Omega, L^{1}(\Omega)) \cap
L^{p}(\Omega, L^{p'}(\Omega))$ for some $1\leq p <\infty$. 


%

        Moreover assume that  the locally 
	Lipschitz function $f$ satisfies that 
	$f(\cdot,0)\in L^{\infty}
	(\Omega)$, and
	\begin{equation}\label{cond_f1}
	\displaystyle\frac{\partial f}{\partial s}(\cdot,s)\le \beta(\cdot)\in L^{\infty}
	(\Omega) , \quad s\in \R, 
	\end{equation}
	and for some $1< \rho<\infty$ 
	\begin{equation}\label{cond_f2}
	\displaystyle\left|\frac{\partial f}{\partial s}(\cdot,s)\right|\le C(1+|s|^
	{\rho-1}) , \quad s\in \R . 
	\end{equation}
	 Then equation (\ref{nonlinear_loc_lip}) with initial data 
	$u_0 \in X= L^{p\rho}(\Omega)$  has a unique global solution  
	given by the Variation of Constants Formula 
	\begin{equation}\label{variation_constant_form_f_loc_lip}
		u(\cdot,t)=e^{Lt}u_0+\displaystyle\int_0^t e^{L(t-s)}f(\cdot,u
		(\cdot,s))\,ds,
	\end{equation}
	with $u\in\mathcal{C}\big([0,T], L^{p\rho}(\Omega)\big)\cap \mathcal{C}^1
	\big([0,T], L^{p}(\Omega)\big)$ for all $T>0$, 
	and it is a strong solution in $L^{p}(\Omega)$.
\end{theorem}
\begin{proof}
Observe that from the assumptions on $J$ we get
$K\in\mathcal{L}(L^{\infty}(\Omega), L^{\infty}(\Omega))$ and 
$K\in\mathcal{L}(L^p(\Omega), L^p(\Omega))$. Then by Riesz--Thorin
interpolation we get  $K\in\mathcal{L}(L^{p\rho}(\Omega),
L^{p\rho}(\Omega))$.

	Now we  prove that $f$ satisfies (\ref{hyp_f_Cs_2_D}). 
	Let $s>0$ be arbitrary. Integrating (\ref{cond_f1}) in $[0,s]$, and 
	multiplying   by $s>0$, we obtain 
	$$
	f(\cdot,s)s  \le \beta(\cdot)s^2+f(\cdot,0)s
		\le Cs^2+D|s| , 
	$$
with an analogous argument for  $s<0$. 
Thus, from Proposition \ref{global_sol_initial_data_bounded_2}  we have the existence and 
	uniqueness of solutions for (\ref{nonlinear_loc_lip}) with initial data  	$u_0\in L^{\infty}(\Omega)$. 
	
Denoting $q=p\rho$, since $L^{\infty}(\Omega)$ is dense in $L^q(\Omega)$, we consider a sequence of initial data
$\{u_0^n\}_{n\in\mathbb{N}}\subset L^{\infty}(\Omega)$ such that
$u^n_0\to u_0$ in $L^{q}(\Omega)$ as $n$ goes to $\infty$.  Thanks to Proposition
\ref{global_sol_initial_data_bounded_2}, we know that the
solution of (\ref{nonlinear_loc_lip}) associated to the
initial data $u_0^n\in L^{\infty}(\Omega)$, satisfies
\[
u^n_t(x,t) =L\,u^n(x,t)+f(x,u^n(x,t)).
\]

	We prove first that $\left\{u^n\right\}_{n\in\mathbb{N}}\subset \mathcal{C}
	([0,\infty),L^q(\Omega))$ is a 
	Cauchy sequence in compact sets of $[0,\infty)$. Since 
	\begin{equation}\label{eq_repuesto_3}
	u^k_t(t)-u_t^j(t)=L(u^k-u^j)(t)+f(\cdot,u^k(t))-f(\cdot,u^j(t)), 
	\end{equation}
	multiplying (\ref{eq_repuesto_3}) by $|u^k-u^j|^{q-2}(u^k-u^j)(t)$, and 
	integrating in $\Omega$, we obtain
	\begin{equation}\label{derivative_cauchy_sequences}
	\begin{array}{ll}
	\displaystyle\frac{1}{p}\frac{d}{d\, t} \|u^k(t)-u^j
          \!\!\!\!\!\!& (t)\|^q_{L^q(\Omega)}  =\displaystyle
                        \int_{\Omega} L(u^k-u^j)(t) |u^k-u^j|^{q-2}(u^k-u^j)(t) \smallskip\\
	& + \displaystyle\int_{\Omega}\left(f(\cdot,u^k(t))-f(\cdot,u^j(t))\right)|u^k\!-u^j|^
	{q-2}(u^k-u^j)(t).
	\end{array}
	\end{equation}
Denoting $w(t)= u^k(t)-u^j(t)$ and $g(w)= |w|^{q-2}w\in
L^{q'}(\Omega)$, 
we write 
	\begin{equation}\label{derivative_cauchy_sequences_001}
		L w(t)=(K -h_0(\cdot)) w(t)+(h_0(\cdot) -h(\cdot)) w(t)
	\end{equation}
and then since  $J(x,y)=J(y,x)$ and $K\in\mathcal{L}(L^q(\Omega), L^q(\Omega))$, we get 
	\begin{displaymath}
\int_{\Omega} (K-h_0I)\,w g(w)\,dx = -\frac{1}{2}\int_{\Omega}\int_{\Omega}J(x,y)(w(y)-w(x))(g(w)
	(y)-g(w)(x))dy\,dx.	
	\end{displaymath}
Since  $J$ is nonnegative and $g(w)= |w|^{q-2}w$ is increasing, then
we obtain  
	\begin{equation}\label{des1}
	\begin{array}{l}
	\!\!\! \displaystyle\int_{\Omega} (K-h_0I)\,w g(w) \,dx \displaystyle =
	-\frac{1}{2}\int_{\Omega}\!\int_{\Omega}J(x,y)(w(y)\!-\!w(x))(g(w)
	(y)-g(w)(x))dy\,dx\le 0.
	\end{array}
	\end{equation}
	Moreover, $h,\, h_0\in L^{\infty}(\Omega)$, then from the second part on the 
	right hand side of 
	(\ref{derivative_cauchy_sequences_001}) we obtain in 
	(\ref{derivative_cauchy_sequences}) 
	\begin{equation}\label{des1_2}
	\displaystyle\int_{\Omega} (h_0(x)-h(x))\,|w|^q(x)\,dx \le C\|w\|^q_{L^q(\Omega)}.
	\end{equation}

	On the other hand, thanks to (\ref{cond_f1}) and the mean 
	value Theorem, there exists $\xi=\xi(x,t)$, 
	such that, the second term on the right hand side of 
	(\ref{derivative_cauchy_sequences}) satisfies that
	\begin{equation}\label{des2}
	\begin{array}{l}
	\!\displaystyle\int_{\Omega}\!\!\left(f(\cdot,u^k(t))\!-\!f(\cdot,u^j(t))\right)|u^k-
	u^j|^{q-2}(u^k-u^j)(t)\displaystyle \!=\!\!\int_{\Omega}{\partial f\over \partial 
	u}(\cdot,\xi) |w|^q\le 
	\|\beta\|_{L^{\infty}(\Omega)}\|w\|^q_{L^q(\Omega)}.
	\end{array}
	\end{equation}
Therefore, thanks to (\ref{derivative_cauchy_sequences}), (\ref{des1}), 
	(\ref{des1_2}),  and (\ref{des2}), we obtain
	\[\frac{d}{d\, t} \|u^k(t)-u^j(t)\|^q_{L^q(\Omega)}\le C \|u^k(t)-u^j(t)\|^q_{L^q
	(\Omega)}\] 
and  Gronwall's inequality gives 
	\begin{equation}\label{cons_of_gronwall}
		\|u^k(t)-u^j(t)\|^q_{L^q(\Omega)}\le e^{Ct}\|u_0^k-u_0^j\|^q_{L^q(\Omega)},
	\end{equation}
	and from this 
	\begin{equation}\label{ref_exis_un_u_0_l_p_1}
		\sup\limits_{t\in[0,T]}\|u^k(t)-u^j(t)\|^q_{L^q(\Omega)}\le C(T)\|
		u_0^k-u_0^j\|^q_{L^q(\Omega)}.
	\end{equation}

Now the  right hand side of  (\ref{ref_exis_un_u_0_l_p_1}) goes to 
	zero as $k$ and $j$ go to $\infty$. 
	Therefore we have that $\{u^n\}_n\subset \mathcal{C}([0,\infty),L^q
	(\Omega))$ is a Cauchy sequence in compact sets of
        $[0,\infty)$, and then  the limit of the sequence $\{u^n\}_n$ in $\mathcal{C}
	([0,T],L^q(\Omega))$ for any $T>0$,
	$$u(t)=\lim\limits_{n\to\infty}u^n(t)$$ 
exists  and it is independent 
	of the sequence $\{u_0^n\}_n$. 

From  (\ref{cond_f2}), and since $|\Omega|<\infty$, we have  that
$f:L^{p\rho}(\Omega)\to L^{p}(\Omega)$ 	is Lipschitz in bounded sets
of $L^{p\rho}(\Omega)$.  
Then 	for any $T>0$,
\begin{equation}\label{cota_f(u_n)}
	f(u^n)\rightarrow f(u)\quad \mbox{ in }
        \mathcal{C}([0,T],L^{p}(\Omega)) . 
	\end{equation}
	Since $L\in\mathcal{L}(L^{p}(\Omega), L^{p}(\Omega))$,  then there exists 
	$\delta>0$, such that 
$\|e^{Lt}\|_{\mathcal{L}(L^p(\Omega))}\le C_0e^{\delta t}$.
	Thus
	\begin{equation}\label{ref_ex_L_p_ref_001}
		\begin{array}{ll}
		\displaystyle\bigg\|\int_0^t\!\!e^{L(t-s)}f(\cdot,u^n(s))ds - 
		\displaystyle\int_0^te^{L(t-s)}\!\!\!\!\!&f(\cdot,u(s))ds\bigg\|_{L^p(\Omega)}
		\displaystyle \le \int_0^t e^{\delta s}\,\|f(\cdot,u^n(s))-f(\cdot,u(s)))\|_{L^p
		(\Omega)}\,ds.
		\end{array}
	\end{equation}
	Taking supremums in $[0,T]$ in (\ref{ref_ex_L_p_ref_001}), and 
	from (\ref{cota_f(u_n)}) we obtain 
	\[\displaystyle \int_0^t\!\!e^{L(t-s)}f(\cdot,u^n(s))ds \rightarrow \!\!\int_0^t\!\!e^{L(t-s)}f
	(\cdot,u(s))ds\;\mbox{ in } \mathcal{C}([0,T],L^p(\Omega)),\,\forall T\!>\!0.\]

Also, since $u^n_0\to u_0$ in $L^{p\rho}(\Omega)$ as $n \to \infty$ we
have $e^{Lt}u_0^n\rightarrow e^{Lt}u_0$ in   $\mathcal{C}\left([0,T], 
	L^{p\rho}(\Omega)\right)$ for all $T>0$
and using $\displaystyle \int_0^te^{L(t-s)}f(\cdot,u^n(s))ds=u^n(t)-e^{Lt}u_0^n $, 
passing to the limit, we get 
\[
\displaystyle \int_0^te^{L(t-s)}f(\cdot,u^n(s))ds\;\rightarrow \int_0^te^{L
	(t-s)}f(\cdot,u(s))ds = u(t)-e^{Lt}u_0
\]
in $\mathcal{C}\left([0,T], L^{p\rho}(\Omega)\right)$ for any $T>0$. 
Hence, $u \in \mathcal{C}([0,T], L^{p\rho}(\Omega))$ satisfies (\ref{variation_constant_form_f_loc_lip}).

Consider now  $g(t)=f(\cdot,\, u(t))$. Since $u:[0,T]\mapsto L^{p\rho}(\Omega)
	$ is continuous, and $f:L^{p\rho}(\Omega)\mapsto L^{p}(\Omega)$ is 
	continuous, we have that $g:[0,T]\mapsto L^{p}(\Omega)$ is continuous. 
	Moreover, $L\in\mathcal{L}(L^{p}(\Omega), L^{p}(\Omega))$, then, thanks 
	to \cite[Th 2.9, p. 109]{Pazy}, we have that 
	$u\in \mathcal{C}^1([0,T], L^{p}	(\Omega))$ and it is a strong solution of
        (\ref{nonlinear_loc_lip})  in $L^{p}(\Omega)$.
	
	Finally, let us prove the uniqueness of  solutions of   (\ref{nonlinear_loc_lip}) with 
	initial data $u_0 \in L^{p\rho}(\Omega)$, such that $u\in\mathcal{C}\left([0,T], 
	L^{p\rho}(\Omega)\right)	\cap \mathcal{C}^1([0,T],
      L^{p}(\Omega))$ for any  $T>0$, is a strong solution of
      (\ref{nonlinear_loc_lip}) and  is given by the variations of
      constants formula (\ref{variation_constant_form_f_loc_lip}).  
Indeed there exist two such solutions $u$ and $v$, following 
	the steps of this proof from (\ref{eq_repuesto_3}) to (\ref{cons_of_gronwall}),
	 replacing $u^k$ for $u$ and $u^j$ for $v$, we obtain 
	\begin{displaymath}
		\|u(t)-v(t)\|^{p\rho}_{L^{p\rho}(\Omega)}\le e^{Ct}\|u(0)-v(0)\|^{p\rho}_{L^{p\rho}(\Omega)}.
	\end{displaymath}
From this, uniqueness follows. 
\end{proof}

In the following Corollaries we enumerate the  monotonicity properties that 
are satisfied for the solution of (\ref{nonlinear_loc_lip})
constructed in Theorem \ref{global_solutions_f_growth}. 

\begin{corollary} {\bf (Weak, strict  and strong comparison principles) }

Under the   assumptions
of Theorem \ref{global_solutions_f_growth},  for any initial data
$u_0\in X= L^{p\rho} (\Omega)$ 

\noindent i) If $J\geq 0$ then  problem (\ref{nonlinear_loc_lip}) satisfies a weak
  comparison principle. 
	
\noindent ii) If, additionally,  $\Omega$ is $R$-connected and $J$ satisfies 
hypothesis (\ref{eq:J_positive}) then problem (\ref{nonlinear_loc_lip}) satisfies a strict 
  comparison principle. 

\noindent iii) Finally, if  moreover  $|\Omega|<\infty$, $\Omega$ is  $R$--connected and  the measure satisfies
  (\ref{eq:measure_non_degenerates}) and $J$ satisfies
  (\ref{eq:J_strict_positive})  then problem (\ref{nonlinear_loc_lip}) satisfies a strong 
  comparison principle.


\end{corollary}
\begin{proof}
	Given $u_0,\,u_1\in L^{p\rho}(\Omega)$, with $u_0\ge u_1$,
        since  $L^{\infty}(\Omega)$
	is dense in $L^{p\rho}(\Omega)$, then we choose two 
	sequences $\{u_0^n\}_{n\in\mathbb{N}}$ and $\{u_1^n\}_{n\in\mathbb{N}}$ in 
	$  L^{\infty}(\Omega)$ that converge to the initial data $u_0$ 
	and $u_1$ respectively, and such that $u_0^n\ge u_1^n$,   for
        all $n\in\mathbb{N}$.
	
	Thanks to Corollary \ref{comparison_f1f2_Loc_Lipschitz} we know 
	that the associated solutions satisfy $u_n^0( t)\ge u_n^1(
        t)$,  for $t\ge 0$ and $n\in\mathbb{N}$.
	From Theorem \ref{global_solutions_f_growth}, we know that $u_n^i(t)$ 
	converges to $u^i(t)$, for $i=0,1$ in $\mathcal{C}([0,T], L^{p\rho}(\Omega))
	$.  Therefore $u^0( t)\ge u^1( t)$,  for $t\ge 0$.
As in Proposition \ref{comparison_G1G2_Lipschitz}  we
arrive to $u^0( t) > u^1( t)$   or $\inf_{\Omega} (u^0( t) - u^1( t)) >0$  for all $t>0$
        respectively. 
\end{proof}

\begin{corollary}
{\bf (Weak and strong parabolic maximum principle) } 

Under the assumptions in Theorem \ref{global_solutions_f_growth}, assume
moreover that $f(\cdot, 0)\ge 0$, Then for initial data in   $X=
L^{p\rho}(\Omega)$

\noindent i) If $J\geq 0$ then  problem (\ref{nonlinear_loc_lip})satisfies a weak
  maximum  principle. 
	
\noindent ii) If, additionally,  $\Omega$ is $R$-connected and $J$ satisfies 
hypothesis (\ref{eq:J_positive}) then problem (\ref{nonlinear_loc_lip})
satisfies a strict 
  maximum  principle. 

\noindent iii) Finally, if  moreover  $|\Omega|<\infty$, $\Omega$ is  $R$--connected and  the measure satisfies
  (\ref{eq:measure_non_degenerates}) and $J$ satisfies
  (\ref{eq:J_strict_positive})  then problem (\ref{nonlinear_loc_lip})
satisfies a strong 
maximum  principle.


\end{corollary}

\begin{corollary}\label{ordered_solutions_supersolutions_f_growth}

Under the assumptions in Theorem \ref{global_solutions_f_growth}, let $u( t,\, u_0)$ be the solution to 
	(\ref{nonlinear_loc_lip}) with initial data $u_0\in
        L^{p\rho}(\Omega)$ and $\bar{u}( t)$ be a supersolution
	 to (\ref{nonlinear_loc_lip}) in $[0, T]$. 

If $\bar{u}(0)\ge u_0$, then 
	$$\bar{u}( t)\ge u( t,\, u_0),\quad \mbox{ for all }\; t\in [0, T].$$
	The same is true for subsolutions with reversed inequality.
\end{corollary}

%
%
%
%
%

\section{Asymptotic estimates}
\label{asymptotic_estimates}

In this section we will show how structure conditions on the nonlinear
term, which correspond to some sign condition at infinity, allow to
obtain suitable estimates on the solutions of
\begin{equation}\label{nonlinear_loc_lip_asymp_est}
	\left\{\!\!
	\begin{array}{rl}
	u_t(x,t)\!\! & =Ku(x,t)\!+\!f(x,u(x,t)),  \quad x\in\Omega,\, t>0,\\
	u(x,0)  \!\!& =u_0(x) , \quad x\in\Omega 
	\end{array}
	\right.
      \end{equation}
      where we assume $|\Omega|<\infty$. Notice that we have set,
      without loss of generality, $h(x)=0$ and we will assume $f(x,s)$ is
locally Lipchitz in $s\in \R$, uniformly in
$x\in \Omega$. If we work in $X=C_{b}(\Omega)$ then we will additionally requiere $f(x,s)$
continuous in $(x,s)\in \Omega \times \R$.


Notice that, from the previous sections,  for the existence of solutions of
(\ref{nonlinear_loc_lip_asymp_est}) 
we will assume either one of the following situations:

\begin{align}  \label{eq:nonlinear_case_1}  
\parbox{0.85\textwidth}{As in Proposition \ref{existence_global_lipschitz_G},
$0\leq J\in L^{p}(\Omega, L^{p'}(\Omega))$ for some $1\leq p \leq \infty$
and then denote $X= L^{p}(\Omega)$, or $J\in C_{b}(\Omega, L^{1}(\Omega))$, 
then denote $X= C_b(\Omega)$.  Also  $g=f(\cdot,0) \in
X$  and  $f(x,s)$ globally Lipschitz in $s$. }
\end{align}

\begin{align} \label{eq:nonlinear_case_2}
\parbox{0.85\textwidth}{As in Proposition
  \ref{global_sol_initial_data_bounded_2},  $0\leq J\in L^{\infty}(\Omega,
  L^{1}(\Omega))$ and then denote $X= L^{\infty}(\Omega)$, or $J\in
C_{b}(\Omega, L^{1}(\Omega))$ and then denote $X= C_b(\Omega)$. Also
  $g=f(\cdot,0) \in X$ and $f(x,s)$ locally  Lipschitz 
  in $s$ and satisfies (\ref{hyp_f_Cs_2_D}).}
\end{align}

\begin{equation}\label{eq:nonlinear_case_3}
\begin{gathered}
\parbox{0.85\textwidth}{As  in Theorem \ref{global_solutions_f_growth},   $J(x,y)=J(y,x)$,  
$0\leq J\in L^{\infty}(\Omega, L^{1}(\Omega)) \cap
L^{p}(\Omega, L^{p'}(\Omega))$ for some $1\leq p <\infty$ and for some
$1< \rho<\infty$}
\\ 
	\displaystyle\frac{\partial f}{\partial s}(\cdot, s)\le \beta(\cdot)\in L^{\infty}
	(\Omega),    \quad 
\displaystyle\left|\frac{\partial f}{\partial s}(\cdot,s)\right|\le
C(1+|s|^{\rho-1}) ,  \\
\parbox{0.85\textwidth}{$g=f(\cdot,0) \in L^{\infty} (\Omega)$, and $X= L^{p\rho}(\Omega)$.}
\end{gathered}
\end{equation}


\noindent Hence, in all three cases the global solutions of
(\ref{nonlinear_loc_lip_asymp_est}) for $u_{0} \in X$ allows us to
define a nonlinear semigroup of solutions  in $X$ by 
\begin{equation} \label{eq:nonlinear_semigroup}
  S(t)u_{0} = u(t, u_{0}), \qquad t \geq 0, \quad u_{0} \in X . 
\end{equation}

\begin{proposition}\label{linear_bound_on_u}


With either one of  the assumptions above (\ref{eq:nonlinear_case_1}),
(\ref{eq:nonlinear_case_2}) or (\ref{eq:nonlinear_case_3}), 
   there exist  $C, D \in L^{\infty} (\Omega)$ such that
   $|g(x)|\leq D(x)$ and 
\begin{equation}\label{structure_condition_on_f}
	f(x,s)s \le C(x)s^2+D(x)|s|,\qquad  s \in \R, \quad x \in
        \Omega 
\end{equation}
and moreover assume   $C, D \in C_{b}(\Omega)$ if $X=C_{b}(\Omega)$.

	Let $\mathcal{U}(t)$ be the solution of 
	\begin{equation}\label{eq_of_supersol}
	\left\{
	\begin{array}{ll}
	\mathcal{U}_t(x,t)=K\,\mathcal{U}(x,t)+C(x)\mathcal{U}(x,t)+D(x)
	 & x\in\Omega,\, t>0, \\
	\mathcal{U}(x,0)=|u_0(x)| & x\in \Omega . 
	\end{array}
	\right.
	\end{equation} 
	Then the solution, $u$, of (\ref{nonlinear_loc_lip_asymp_est}), satisfies that 
        \begin{equation}
          \label{eq:u_below_U}
	|u(t)|\le \mathcal{U}(t),\;\mbox{ for all }\; t\ge 0.          
        \end{equation}

\end{proposition}
\begin{proof}
Observe that in  case of assumptions \eqref{eq:nonlinear_case_1} above, 
(\ref{structure_condition_on_f}) is  satisfied with $C(x)= L_{0}$ where $L_{0}$ is  the Lipschitz constant of  $f(x,s)$ in
$s$, and $D(x) = |g(x)|$.  On the other hand, for assumptions
(\ref{eq:nonlinear_case_2}), we have \eqref{structure_condition_on_f} since $f$ satisfies 
(\ref{hyp_f_Cs_2_D}). Finally in case of assumption (\ref{eq:nonlinear_case_3}) then
(\ref{structure_condition_on_f}) is  satisfied with $C(x)=\beta(x)$ and $D(x) = |g(x)|$.

	Now we prove that  the solution of
        (\ref{eq_of_supersol}) is nonnegative.  In fact we  know that,
        denoting $L_C=K + CI$, 
  	\begin{equation}\label{VCF_for_U}
	\mathcal{U}(t)=e^{L_Ct}|u_0|+\displaystyle\int_0^t
        e^{L_C(t-s)}D  \, ds . 
	\end{equation}
	Since $|u_0|, D\ge 0$,  
then we can apply Proposition
\ref{positive_solution_with_u_0_positive} and then 
	we have that $e^{L_Ct}|u_0|\ge 0$ and $	e^{L_C(t-s)}D\ge 0$
        for $t\geq 0$. 
	 Thus, we have that $\mathcal{U}(t)$ is nonnegative for all $t\ge 0$. 

	Now, we prove that $\mathcal{U}$ is a supersolution of
        (\ref{nonlinear_loc_lip_asymp_est}). First, since in any of
        the cases (\ref{eq:nonlinear_case_1}),
        (\ref{eq:nonlinear_case_2}) or (\ref{eq:nonlinear_case_3})  we have   $D\in X$
        then $\mathcal{U} \in \mathcal{C}([0,\infty), X)$. 
Now since $\mathcal{U}$ is nonnegative and $f$ 
	satisfies (\ref{structure_condition_on_f}),  we obtain
	\[
	K\mathcal{U} +f(\cdot,\mathcal{U}) \le K\mathcal{U} +C \mathcal{U}+
	D =\mathcal{U}_t.
	\]
	Moreover $u_0\le|u_0|=\mathcal{U}(0)$, then  from either
        Proposition \ref{ordered_solutions_supersolutions_G_Lipschitz} or  Corollary 
	\ref{ordered_solutions_supersolutions_f_Loc_Lipschitz}
        or   \ref{ordered_solutions_supersolutions_f_growth} 
	we have $ u(t)\le \mathcal{U}(t)$ for $t\geq
        0$. 
	Arguing analogously for $-\mathcal{U}(t)$ we obtain
        $-\mathcal{U}(t)\le u(t)\le \mathcal{U}(t)$ for $t\geq 0$ and
        thus the result. 
\end{proof}
%
%

%

Now, we obtain  asymptotic estimates on  the solutions of
(\ref{nonlinear_loc_lip_asymp_est}). 
\begin{proposition}\label{asymp_estimate_phi_nonnegative}
	Let  $X, J, h$  and 
        $f:\Omega\times\R\to\R$ be as in Proposition
        \ref{linear_bound_on_u}. If $X=L^{p}(\Omega)$
        with $1\leq p<\infty$ we furthermore assume that $J\in BUC(\Omega,
        L^{p'}(\Omega))$ whence 
        $K\in\mathcal{L}(L^{p}(\Omega),L^{\infty}(\Omega))$ is
        compact. 
%
%
        Finally assume $f$
        satisfies         (\ref{structure_condition_on_f})  and 
	\begin{equation}\label{semigroup_LC_decays}
	\sup Re\Big( \sigma_{X}(K + CI ) \Big)  \leq  -\delta<0 . 
	\end{equation}

	Then there exists a unique	 equilibrium solution, $\Phi
        \in X$, associated to  
	(\ref{eq_of_supersol}), that is, a solution of 
	\begin{equation}\label{equation_PHI}
		K\Phi + C(x)\Phi + D(x) =0, \quad x\in \Omega 
	\end{equation}  
which moreover satisfies  $\Phi\in L^{\infty}(\Omega)$ and $\Phi\ge
0$. 
If additionally $\Omega$ is  $R$--connected and  the measure satisfies
  (\ref{eq:measure_non_degenerates}) and $J$ satisfies
  (\ref{eq:J_strict_positive}) then, either $D=0$ and then $\Phi=0$ or
  $\inf_{\Omega}\Phi  >0$ in $\Omega$.

Also,  for any  $u_0\in X$,   the solution $u$ of 
	(\ref{nonlinear_loc_lip_asymp_est})  satisfies that
	\[ 
\limsup_{t\to\infty} \|u(t,u_0)\|_{X}\le
          \|\Phi\|_{X}
        \]
        and
        \begin{displaymath}
          \Big(|u(t)|-\Phi \Big)_{+} \to 0 \quad \mbox{in $X$ as $t\to
          \infty$.}
        \end{displaymath}

If moreover  $u_{0} \in L^{\infty}(\Omega)$ then 
\begin{equation} \label{eq:uniform_limit}
  \limsup_{t \to \infty} |u(x, t, u_{0})| \leq \Phi(x), \quad
  \mbox{uniformly in $x\in \Omega$} 
\end{equation}
and uniformly in $u_{0}$ in a bounded set of $L^{\infty}(\Omega)$. 
In particular, any equilibrium, that is, any constant in time solution of
(\ref{nonlinear_loc_lip_asymp_est}), $\varphi \in X$, satisfies
$\varphi  \in L^{\infty}(\Omega)$ and 
\begin{displaymath}
  |\varphi(x)| \leq \Phi(x), \quad a.e. \ x\in \Omega .
\end{displaymath}

Finally if $|u_{0}| \leq \Phi$ then $  |u(t,u_{0})| \leq \Phi$ for
$t\geq 0$. 

\end{proposition}
\begin{proof}
	First of all observe that if   $X=L^{p}(\Omega)$
        with $1\leq p<\infty$, the additional assumption  $J\in BUC(\Omega,
        L^{p'}(\Omega))$ and  Proposition
        \ref{independent_spectrum_K_h} imply that the spectrum
        $\sigma_{X}(K+CI )$ coincides with the spectrum in
        $L^{\infty}(\Omega)$.

	Now from  (\ref{semigroup_LC_decays}), we have 
	that 
	$0$ does not belong to the spectrum of $L_C=K+CI $ hence, it is 
	invertible. 
	Thus, the solution $\Phi$ of (\ref{equation_PHI}) in $X$  is unique. 
	On the other hand, since $D\in L^{\infty}(\Omega)$ and
        $L_C$ is linear and continuous and invertible in $L^{\infty}(\Omega)$,  from equation
        (\ref{equation_PHI}) we also get $\Phi\in L^{\infty} 
	(\Omega)$.
	
	Now, we prove that $\Phi$ is nonnegative. Observe that since
        we know that $\mathcal{U}(t) \geq 0$ satisfies
        (\ref{VCF_for_U}) and,
        thanks  to (\ref{semigroup_LC_decays}),  we have 
	that $\|e^{L_Ct}\|_{\mathcal{L}(X, X)}\le  Me^{-\delta t}$, then
        the limit $\lim\limits_{t\to\infty}\mathcal{U}(t)=\displaystyle\int_0^{\infty} e^{L_
		{C}s}D\,ds \geq 0$ 
exists in $X$, which, from  Lemma \ref{limit_es_equilibrio} below, is an equilibrium of
        (\ref{eq_of_supersol}). Since  this problem has a unique
        equilibrium we obtain  $\Phi = \int_0^{\infty}
        e^{L_{C}s}D\,ds\geq 0$.  

If additionally $\Omega$ is  $R$--connected and  the measure satisfies
  (\ref{eq:measure_non_degenerates}) and $J$ satisfies
  (\ref{eq:J_strict_positive}) then  either $D=0$ and then $\Phi=0$
        or by the strong maximum principle in Theorem
        \ref{thr:Lambda},  $\inf_{\Omega}\Phi        >0$ in $\Omega$.

Now, since  (\ref{eq_of_supersol}) is a linear non--homogenous problem
we can  write
        \begin{equation}
          \label{eq:decomposition_U}
|u(t)|\leq         \mathcal{U}(t)=\Phi+e^{L_Ct}(|u_0|-\Phi)
        \end{equation}
 for any $u_{0} \in   X$. 
From this, 	we obtain 
	\begin{equation}\label{cota_exp_asymp_est_1}
	\|u(t)\|_{X}\le \|\mathcal{U}(t)\|_{X}
	\le \|\Phi\|_{X}+ Me^{-\delta t}\|(|u_0|-\Phi)\|_{X} . 
	\end{equation}
	Since $\delta>0$, then we have
        $\limsup_{t\to\infty} \|u(t)\|_{X}\le  \|\Phi\|_{X}$. 
        Also (\ref{eq:decomposition_U}) gives    
        \begin{displaymath}
       \Big( |u(t,u_{0})| -\Phi\Big)_{+}\leq \Big(e^{L_C
         t}(|u_0|-\Phi)\Big)_{+}  \to 0 \quad
        \mbox{in $X$ as  $t\to \infty$. }   
        \end{displaymath}

        On the other hand, from (\ref{eq:u_below_U}), if $u_{0}\in
L^{\infty}(\Omega)$ we have $\mathcal{U}(t) \to \Phi$ in
$L^{\infty}(\Omega)$ and we get (\ref{eq:uniform_limit}). The result
for the equilibria is then immediate. 

Finally, if $|u_{0}| \leq \Phi$ we get $|u(t,u_{0})|\le
\Phi+e^{L_C t}(|u_0|-\Phi) \leq \Phi$ for $t\geq 0$. 
\end{proof}

\begin{remark}




Notice that for  a given nonlinear term there might be however many
different choices of $C,D$ satisfying
(\ref{structure_condition_on_f}). See Section \ref{logistic_nonlin} for the case of a logistic type nonlinearity.
\end{remark}

\begin{remark}

Observe that, using the notations in Theorem \ref{thr:Lambda},
condition (\ref{semigroup_LC_decays}) is equivalent to $\Lambda(C):=\Lambda(K+CI)
<0$ and then
Proposition \ref{prop:citeria_4_sign_Lambda} and Corollary
\ref{cor:theshold_h0} provide sufficient conditions for that.

\end{remark}
%


Now we prove the lemma used above. 

\begin{lemma}\label{limit_es_equilibrio}
	Let $X$ be a Banach space, and let $S(t):X\to X$ be a 
	continuous semigroup. Assume that $u_0,\,v\in X$ satisfy that 
	$S(t)u_0\to v$ in $X$ as $t\to\infty$. Then $v$ is an equilibrium point 
	for $S(t)$.
\end{lemma}     
\begin{proof}
	Since $v=\lim\limits_{t\to\infty}S(t)u_0$. Then applying $S(s)$ for 
	$s>0$, and using the continuity of $S(t)$ for $t>0$,
        $S(s)v=S(s)\lim\limits_{t\to\infty}S(t)u_0=\lim\limits_{t\to\infty}S(s+t)
        u_0=v$.  
	Then $v$ is an equilibrium point.
\end{proof}

The next result gives some information about the asymptotic behavior
of solutions.

\begin{corollary} \label{cor:pointwise_asymptotic_estimate}
  Under the assumptions of Proposition
  \ref{asymp_estimate_phi_nonnegative},   for any  $u_0\in X$  and  any  sequence $t_{n} \to \infty$ there exists a subsequence
(that we denote the same) such that $u(t_{n}, u_{0})$ converges 
weakly in $X$ (or weak* if $X=L^{\infty}(\Omega)$) to some  bounded function
$\xi$ such that  $| \xi(x) |\leq \Phi(x)$ a.e.  $x\in \Omega$ and 
\begin{equation} \label{eq:pointwise_limit}
\limsup_{n} |u(t_{n}, x, u_{0})| \leq \Phi(x), \quad a.e. \ x\in
\Omega.
\end{equation}

\end{corollary}
\begin{proof}
 From (\ref{cota_exp_asymp_est_1}) we have that 
 $\{u(t,u_{0}), \ t\geq 0\}$ is bounded in $X$. Therefore, taking
 subsequences if necessary, we can assume that $u(t_{n},u_{0})$ converges weakly to $\xi
\in X$ and  $e^{L_Ct_{n}}(|u_0|-\Phi)\to 0$ a.e. in $\Omega$. Hence
from (\ref{eq:decomposition_U}) we get the
result.
\end{proof}

\section{Extremal equilibria}\label{section_extr_equil}

In this section we prove that (\ref{nonlinear_loc_lip_asymp_est}) has two ordered
extremal equilibria $\varphi_m\le \varphi_M$, stable from below and
from above, respectively, that enclose the asymptotic behavior of all
solutions.


Observe that from  Proposition
\ref{asymp_estimate_phi_nonnegative} and Corollary
\ref{cor:pointwise_asymptotic_estimate}, in order to analyze the 
asymptotic behavior of solutions,  one can always assume that
\begin{displaymath}
|u_{0}|\leq \Phi \in L^{\infty}(\Omega), \quad \mbox{and then} \quad   |u(t,u_{0})|\leq
\Phi \in L^{\infty}(\Omega) \quad \mbox{for all $t\geq 0$} . 
\end{displaymath}
In particular we can always assume $f(x,s)$ is globally Lipschitz in
its second variable and the nonlinear semigroup in
(\ref{eq:nonlinear_semigroup}) is continuous in the
norm of $L^{q}(\Omega)$ for any $1\leq q \leq \infty$.

\begin{theorem}\label{two_extremal_equilibria}

Under the assumptions of Proposition
\ref{asymp_estimate_phi_nonnegative},  there exist  two 
	ordered bounded extremal equilibria  of  the problem 
	(\ref{nonlinear_loc_lip_asymp_est}), $\varphi_m\le 
	\varphi_M $, with $|\varphi_{m}|, |\varphi_{M}| \leq \Phi$, 
	such that any other equilibria $\psi$ of (\ref{nonlinear_loc_lip_asymp_est}) 
	satisfies $\varphi_m\le \psi \le \varphi_M$. 
	
	Furthermore, the set 
	$\left\{v\in L^{\infty}(\Omega):\varphi_m\le v\le \varphi_M\right\}$ 
	attracts the dynamics of the solutions $u(t,u_0)$ of the problem 
	(\ref{nonlinear_loc_lip_asymp_est}), in the sense that for
        each  $u_0\in L^{\infty}	(\Omega)$, there exist $\underline{u}(t)$ and $\overline{u}(t)$ 
	in $L^{\infty}(\Omega)$ such that $\underline{u}(t)\le u(t,u_0)\le 
	\overline{u}(t)$, and 
	\[
	\begin{array}{l}
		\lim\limits_{t\to\infty}\underline{u}(t)=\varphi_m , \qquad\qquad
		\lim\limits_{t\to\infty}\overline{u}(t)=\varphi_M
	\end{array}
	\]
	in $L^{q}(\Omega)$ for any  $1\le q< \infty$. 
\end{theorem}
\begin{proof}
  From  (\ref{eq:decomposition_U}), 
	since $\|e^{L_Ct}\|_{\mathcal{L}(L^{\infty}(\Omega))}\le M e^{-\delta t}$, 
	with $\delta>0$, and  $u_0\in L^{\infty}(\Omega)$ then 
	for  $\varepsilon>0$ there exists $T(u_0)>0$ such that $\|e^{L_Ct}(|u_0|-\Phi)\|_{L^{\infty}(\Omega)}
		<\varepsilon$, for $ t\ge T(u_0)$. 
	Then  (\ref{eq:decomposition_U})  gives $-\Phi-\varepsilon\le u( t,
        u_0)\le \Phi+\varepsilon$ for  $t\ge 	T(u_0)$ 
 which, writing $T=T(u_0)$, we recast  as 
	\begin{equation}\label{eq_equil_3}
		-\Phi-\varepsilon\le S(t+T) u_0 \le \Phi+\varepsilon,\quad \forall t
		\ge0.
	\end{equation}

In particular, for  the initial data $u_0=\Phi+
	\varepsilon$, then  there exists $T=T(\Phi+\varepsilon)$ such that 
	\begin{equation}\label{eq_equil_4}
		-\Phi-\varepsilon\le S(t+T)(\Phi+\varepsilon)\le \Phi+\varepsilon,
		\quad \forall t\ge0.
	\end{equation}
	Now, from the comparison principles in Section
        \ref{exist_uniq_pos_nonlinear} 
	and applying $S(T)$ to (\ref{eq_equil_4}) with $t=0$, we
        obtain 
	 \begin{displaymath} 
		-\Phi-\varepsilon\le S(2T)(\Phi+\varepsilon)\le S(T)(\Phi+
		\varepsilon)\le \Phi+\varepsilon.
	\end{displaymath}
	Iterating this process, we obtain that
	\begin{displaymath} 
		-\Phi-\varepsilon\le S(nT)(\Phi+\varepsilon)\le S((n-1)T)(\Phi+
		\varepsilon)\le \dots \le S(T)(\Phi+\varepsilon)\le \Phi+\varepsilon,\; 
	\end{displaymath}
	for all $ n\in\mathbb{N}$. Thus,
        $\left\{S(nT)(\Phi+\varepsilon)\right\}_{n\in\mathbb{N}} $ is
        a monotonically decreasing sequence bounded from below. From
        the Monotone Convergence Theorem, the sequence converges
        pointwise and in $L^q(\Omega)$, for any $1\le q<\infty$, to
        some function $\varphi_M$, i.e.
	\begin{equation}\label{para_remark_despues_Dini}
	S(nT)(\Phi+\varepsilon)\to \varphi_M\;\mbox{ as }\;n\to 	
	\infty\;\mbox{ in }\; L^q(\Omega). 
	\end{equation}
        From (\ref{eq:pointwise_limit}) we get
        $|\varphi_{M}(x)|\leq \Phi(x)$ in $\Omega$ and 	$\varphi_M\in L^{\infty}(\Omega)$. 



	Now we prove that, in fact, the whole solution $S(t)(\Phi+\varepsilon)$ 
	converges  in $L^q(\Omega)$  to $\varphi_M$ as $t\to\infty$. 
Let $\{t_n\}_{n\in\mathbb{N}}$ be a time sequence tending to
infinity. We can write $t_{n} = k_{n} T + s_{n}$ with integers $k_{n}\to
\infty$ and $0\leq s_{n} < T$. Then from (\ref{eq_equil_4}) we get $
S(s_{n}+T)(\Phi+\varepsilon)\le \Phi+\varepsilon$ 
and then applying $S((k_{n}-1)T)$ to both sides we get 
\begin{equation}\label{ref_1_monotone_conv_2}
 S(t_{n})(\Phi+\varepsilon)\le
S((k_{n}-1)T)(\Phi+\varepsilon) . 
\end{equation}
On the other side, again from (\ref{eq_equil_4}) we also get $
S(2T-s_{n})(\Phi+\varepsilon)\le \Phi+\varepsilon $ 
and then applying $S(t_{n})$ to both sides we get 
\begin{equation}\label{ref_2_monotone_conv_2}
 S((k_{n}+2)T)(\Phi+\varepsilon)\le
S(t_{n})(\Phi+\varepsilon) . 
\end{equation}
        Then, using (\ref{para_remark_despues_Dini}) and 
	taking limits as $n$ goes to infinity in (\ref{ref_1_monotone_conv_2}) 
	and (\ref{ref_2_monotone_conv_2}), 
	we obtain that 
$\lim\limits_{n\to\infty}S(t_n)(\Phi+\varepsilon)=\varphi_M$ in $L^q(\Omega)$ 
        for any $1\leq q<\infty$. Since the previous argument is valid for any sequence 
	$\{t_n\}_{n\in\mathbb{N}}$ we actually have
	\begin{equation}\label{s_t_limit_u_sup}
		\lim\limits_{t\to\infty}S(t)(\Phi+\varepsilon)=\varphi_M\;\;\mbox
		{  in }\; L^q(\Omega) 
	\end{equation}
for any $1\leq q<\infty$. From Lemma \ref{limit_es_equilibrio}
$\varphi_{M}$ is an equilibria. 

Analogously, we obtain the equilibria $\varphi_{m}$ as
$\lim\limits_{t\to\infty}S(t)(-\Phi-\varepsilon)=\varphi_m$ in  $ L^q(\Omega) $
for any $1\leq q<\infty$ and $\varphi_{m} \leq \varphi_{M}$. 

Now,  for a general 
	initial data   $u_0\in L^{\infty}(\Omega)$, from  
	 (\ref{eq_equil_3}), for $T=T(u_0)$
	\begin{equation}\label{ineq_equil_u0_1}
u(t+T,u_{0})=  S(t+T) u_0 \le \Phi+\varepsilon,\quad \forall t
		\ge0.
	\end{equation}
	Letting the semigroup act at time $t$ in (\ref{ineq_equil_u0_1}), we 
	have 
	\begin{equation}\label{eq_equil_7}
u(2t+T,u_{0})=	S(2t+T)u_0\le S(t)(\Phi+\varepsilon):=\overline{u}(t),\qquad \forall t
		\ge 0.
	\end{equation}
	Thanks to (\ref{s_t_limit_u_sup}) we obtain the result. 
	
	Finally, let $\psi \in L^{\infty}(\Omega)$ be another
        equilibrium. From (\ref{s_t_limit_u_sup}) and (\ref{eq_equil_7}) with
        $u_0=\psi$, we have $\psi\le \varphi_M$.  Thus $\varphi_M$ is
        maximal in the set of equilibrium points. 
The results  for $\varphi_m$ can be obtained in an analogous way.
\end{proof}

In particular, we  obtain the one side stability of the extremal equilibria. 

\begin{corollary}\label{extrem_equil_stable_from_above}
	Under the hypotheses of  Theorem 
	\ref{two_extremal_equilibria}, if $u_0\in L^{\infty}(\Omega)$, and $u_0\ge \varphi_M$,
	then 
	$$\lim\limits_{t\to\infty}u(t,u_0) = \varphi_M,$$
	in $L^q(\Omega)$, for any  $1\le q<\infty$, i.e., $\varphi_M$ is ``stable 
	from above''. 
	
	If $u_0\in L^{\infty}(\Omega)$, and  $u_0\le \varphi_m$, then 
	$$
\lim\limits_{t\to\infty} u(t, u_0) = \varphi_m,  
$$ 
in $L^q(\Omega)$, for any  $1\le q<\infty$, i.e.  $\varphi_m$ is ``stable from below''.
\end{corollary}
\begin{remark}
	 If the extremal equilibria was more regular $\varphi_M\in\mathcal{C}_b(\Omega)$, then 
	 the result of the previous Theorem \ref{two_extremal_equilibria} 
	 could be  
	 improved because we would obtain the asymptotic dynamics of 
	 the solution of (\ref{nonlinear_loc_lip_asymp_est}) enter between the 
	 extremal equilibria uniformly on compact sets of $\Omega$. 
	In fact, thanks to  Dini's Criterium (cf.  
	\cite[p. 194]{Bartle}), we have that the monotonic sequence 
	 $S(nT)(\Phi+ \varepsilon)$ in	 (\ref{para_remark_despues_Dini}),  
	 converges 
	uniformly in compact subsets of 
	$\Omega$ to $\varphi_M$ as $n$ goes to infinity and from this,
        as in (\ref{s_t_limit_u_sup}), 
	\begin{displaymath} 
		\lim\limits_{t\to\infty}S(t)(\Phi+\varepsilon)=\varphi_M \;\;\mbox{ in }
		\; L^{\infty}_{loc}(\Omega).
	\end{displaymath}
Then  Theorem \ref{two_extremal_equilibria} and Corollary
\ref{extrem_equil_stable_from_above} can be stated with this
convergence. 
However, since  there is no regularization for the semigroup $S(t)$ associated 
	to (\ref{nonlinear_loc_lip_asymp_est}), we can not assure that 
	$\varphi_M\in\mathcal{C}_b(\Omega)$, as happens for the local 
	reaction diffusion equations. Indeed there are explicit 
        examples of (non isolated and) discontinuous equilibria. See the example below.    
\end{remark}

\begin{example}{\bf (Example of non-isolated and discontinuous  equilibria)}
\label{example_discontinuous_equilibria} 	
Choose $J(x,y)=1$, for all $x, y\in\Omega$. Then  the equilibria
of the problem 
\begin{displaymath}
  u_t(x,t)= \int_{\Omega} (u(y)-u(x)) \, dy+f(u(x)),
\end{displaymath}
satisfy that  
\begin{equation} \label{ex_equil_nonisolated}
   \int_{\Omega} u(y)\, dy = |\Omega|u(x)  -f(u(x)) \quad x\in \Omega. 
\end{equation}
We construct piecewise constant  solutions of (\ref{ex_equil_nonisolated}) in
the following way. Define $g(u)=|\Omega|u- f(u) $. 
Then choose $A \in \R$ and  consider the set of (real)  solutions  of
$g(u)=A$, that is $f(u)= |\Omega|u -A$. Then any piecewise constant function $u(x)$ with values in
the set of solutions of $g(u)=A$ will be an equilibria of the problem,
provided 
 $$
 \int_{\Omega}u(y) \, dy=A . 
$$
For example we can consider  $f\in C^2_b(\R)$ such that it coincides in an 
interval of the form $u\in [-M,M]$, with $M$ large, with the function 
$f_{0}(u)=\lambda u(1-u^2)$, with $\lambda \in
\R$. 
Thus $f$  gives  three constant equilibria $u^{0}_{1} =-1$, $u^{0}_{2}
=0$ and $u^{0}_{3} =1$. 
For this choice of $f$, we have for $u\in [-M,M]$,
$g(u)=
\big(|\Omega| -\lambda \big) u +\lambda u^{3}$. Depending on the
choice of $\lambda$ and $A$ we can assume that there are three
different roots of $g(u)=A$, which we
denote by  $u_1, u_2, u_3 \in  [-M,M]$.
If we divide the set $\Omega$ in three arbitrary subsets
 $\Omega_1,\,\Omega_2,\, \Omega_3$, then we can construct the
piecewise constant  equilibria
$$
\bar u(x)=u_1  \,\chi_{\Omega_1}(x)+
 u_2 \,\chi_{\Omega_2}(x)+u_3 \,\chi_{\Omega_3}(x), \quad x\in \Omega 
$$
provided $u_1  |\Omega_{1}| + u_2 |\Omega_{2}|  + u_3 |\Omega_{3}| =
  A$,  $ |\Omega_{1}| +  |\Omega_{2}|  +  |\Omega_{3}| =
  |\Omega|$. 
This family of equilibria is not isolated in $L^p(\Omega)$, $1\leq
p<\infty$, because by slightly 
changing the  partition (without changing the measure of each set) of
$\Omega$ to new sets denoted by 
 $\widetilde{\Omega}_1$, $\widetilde{\Omega}_2$, 
 and $\widetilde{\Omega}_3$, the piecewise constant  equilibrium
 $\tilde{u}(x)=u_1 \,\chi_{\widetilde{\Omega}_1}(x)+ u_2
 \,\chi_{\widetilde{\Omega}_2}(x)+u_3
 \,\chi_{\widetilde{\Omega}_3}(x) $,
would be  as close as we want, in
 $L^p(\Omega)$, $1\leq p<\infty$,  to the equilibrium $\bar u$. Note
 however that in $L^{\infty}(\Omega)$ the equilibria $\tilde u$ and
 $\bar u$ are not close. 
Finally  note that this construction only imposes
 restrictions on the measures of the sets  $\Omega_1,\,\Omega_2,\,
 \Omega_3$. Thus, once the measures are fixed there are infinitely
 many possibilities to distribute these sets in $\Omega$. 

Notice that the construction above, by slightly changing the sets
$\Omega_1,\,\Omega_2,\, \Omega_3$, shows that different equlibria can
coincide on sets of positive measure. This can not happen in local
reaction diffusion problems like  (\ref{eq:local_RD}) due to the
maximum principle. 
 
Finally, by  shifting the function $f_{0}(u)$ above to the right, e.g. taking
$f_{0}(u) = \lambda (u-3) (1-(u-3)^{2})$, and choosing $\lambda$ and
$A$ properly we can achieve that the three roots of  $g(u)=A$ lie in
$[-M,M]$ and are now positive and so are the piecewise constant
equilibria constructed above. 
\end{example}

We also get the following result that improves Corollary
\ref{cor:pointwise_asymptotic_estimate}.

\begin{corollary}

  Under the assumptions and notations in  Corollary
\ref{cor:pointwise_asymptotic_estimate} and Theorem \ref{two_extremal_equilibria} for
  any $u_{0} \in X$, 
  and any sequence $t_{n} \to \infty$ there exists a
  subsequence (that we denote the same) such that $u(t_{n},u_{0})$
  converges weakly in $X$ (or weak* if $X=L^{\infty}(\Omega)$) to some
bounded   function $\xi$ with $|\xi|\leq \Phi$  and 
\begin{displaymath} 
		\varphi_m(x)\le \xi(x)\le \varphi_M(x)\quad 
		\mbox{for a.e. } x\in \Omega 
\end{displaymath}
\begin{displaymath}
\varphi_{m}(x) \leq \liminf_{n}   u(t_{n}, x, u_{0}) \leq
\limsup_{n} u(t_{n}, x, u_{0}) \leq \varphi_{M}(x), \quad a.e. \ x\in
\Omega. 
\end{displaymath}

\end{corollary}
\begin{proof}
  From  Proposition \ref{asymp_estimate_phi_nonnegative}, 
	we 	know that for any  initial data $u_0\in X$
	\begin{equation}\label{r_1}
		S(t)u_0=u(t,u_0)\le \Phi+e^{L_C t}(|u_0|-\Phi).
	\end{equation}
	Applying the nonlinear semigroup $S(s)$  to (\ref{r_1}), we 
	have that   
	 \begin{equation}\label{r_2}
			S(s)u(t,u_0)=u(t+s,u_0)\le S(s)(\Phi+e^{L_C t}(|u_0|-\Phi)).
	\end{equation}
	Since the semigroup is continuous in $X$ with  respect to 
	the initial data, the right hand side of (\ref{r_2})
        converges  in $X$, as $t\to \infty$ to $S(s)\Phi$.


Then for any sequence $\{t_n\}_{n\in\mathbb{N}} \to \infty$ we can
assume   the  weak limit in 	$X$ of
$\{u(t_n+s,u_0)\}_{n\in\mathbb{N}}$ exists (or weak* if
$X=L^{\infty}(\Omega)$), and we get from (\ref{r_2}) 
%
$\xi(x)\le S(s)\Phi(x)$, for a.e. $x\in\Omega$. 
        From this and  Corollary \ref{extrem_equil_stable_from_above} we
get $\xi(x)\le \lim\limits_{s\to\infty}
		S(s)\Phi(x)= \varphi_M(x)$, a.e. $x\in\Omega$. 
The result for  the minimal equilibrium $\varphi_m$ is analogous. 
\end{proof}

\subsection{Nonnegative solutions}
\label{sec:nonnegative-solutions}

Now we pay attention to solutions with nonnegative initial data.

\begin{proposition}\label{varphi_M_is_positive}

Under the assumptions in Theorem \ref{two_extremal_equilibria}, if 
$$
g(x)= f(x,0)\ge 0, \quad x\in \Omega
$$
then  $\varphi_M\ge 0$ and there exists a minimal nonnegative
equilibrium $0\leq \varphi_{m}^{+} \leq  \varphi_{M}$. Also for any nonnegative nontrivial $u_{0}\geq 0$
\begin{displaymath}
\liminf_{t\to \infty}  u(x, t, u_{0}) \geq  \varphi_{m}^{+} (x), \quad
a.e. \quad x\in \Omega .
\end{displaymath}
Moreover $\varphi_{m}^{+}$ is nontrivial iff   $g$ is not identically
zero and in such a case $\varphi_{m}^{+}$ is stable from below for
nonnegative initial data (see 
Corollary \ref{extrem_equil_stable_from_above}).  
	 
	In particular,   if  $\Omega$ is $R$-connected and $J$ satisfies 
	hypothesis (\ref{eq:J_positive})   
	 then any nontrivial nonnegative equilibria  of 
	 (\ref{nonlinear_loc_lip_asymp_est})  is 
	 in fact strictly positive.
If moreover 
the measure satisfies
  (\ref{eq:measure_non_degenerates}) and $J$ satisfies
  (\ref{eq:J_strict_positive}) then for each nontrivial nonnegative equilibrium
  $\psi$ 
  \begin{displaymath}
    \inf_{\Omega}\psi   >0 .
  \end{displaymath}

Finally, if $g=0$ then the extremal equilibria in Theorem
\ref{two_extremal_equilibria} satisfy $\varphi_{m} \leq 0 \leq \varphi_{M}$. 

\end{proposition}
\begin{proof}
	Since $g \ge 0$, then $0$ is a subsolution of
        (\ref{nonlinear_loc_lip_asymp_est})  for nonnegative initial
        data.  Then using  Corollary
        \ref{ordered_solutions_supersolutions_f_Loc_Lipschitz}   we know that if $u_0\ge 0$, then 
	\begin{equation}\label{u_pos_pq_0_subsol}
		0\le u(x, t, u_0),\quad  x\in\Omega,\,\, t\ge 0.
	\end{equation}
        Moreover if $g = 0$ then $u(t,0)=0=\varphi_{m}^{+}$ is the
        minimal nonnegative equilibrium. On the other hand, if
        $g\neq 0$, then $0\leq \varphi_{M}$ implies $0\leq u(t,0)\leq
        \varphi_{M}$ and is increasing, nonnegative and 
        bounded above and then converges pointwise and in
        $L^{q}(\Omega)$ for any $1\leq q < \infty$ to a positive
        equilibria $\varphi_{m}^{+}$ which is clearly the minimal
        nonnegative equilibrium and is stable from below.

From  Proposition \ref{asymp_estimate_phi_nonnegative} 
	we know that the solution of (\ref{equation_PHI})
        satisfies $0\leq  \Phi\in L^{\infty}(\Omega)$, then from
        (\ref{u_pos_pq_0_subsol}), we have $u(\cdot, t, \Phi)\ge 0$
        for $t\geq 0$. 
	 From Corollary \ref{extrem_equil_stable_from_above},
         $\lim\limits_{t\to\infty} u(\cdot, t, \Phi)=\varphi_M$ in  $L^{q}(\Omega)$ 
for any $1\leq q <\infty$ and then $0\leq \varphi_{m}^{+} \leq \varphi_M$. 
	
	Moreover, if $\psi$ is a nontrivial nonnegative equilibria of 
	(\ref{nonlinear_loc_lip_asymp_est}) 
	and if $J$ satisfies (\ref{eq:J_positive}), then 
	thanks to Corollary \ref{nonlinear_max_pple_f_Loc_Lipschitz} 
	then for  $ t>0$ we have  $\psi=u(\cdot,t,\psi)>0$ or $\inf_{\Omega}\psi        >0$ respectively. 
\end{proof}

Now we consider the case $g(x)=0$ and prove that if $u=0$ is
linearly unstable then there exists a minimal \emph{positive}
equilibrium which is stable from below. 

\begin{proposition} \label{prop:minimal_positive_0_unstable}

 Under the hypotheses of  Theorem \ref{two_extremal_equilibria},
 assume additionaly  
 that   $\Omega$ is
 $R$--connected,  the measure satisfies 
  (\ref{eq:measure_non_degenerates}) and  $J \in BUC(\Omega,
L^{p'}(\Omega))$ satisfies
  (\ref{eq:J_strict_positive}).  

Assume  $g =0$ and that for some $s_{0}>0$ we have  
\begin{displaymath}
  f(x,s) \geq M(x) s , \quad x\in \Omega, \quad 0\leq s \leq s_{0}
\end{displaymath}
with $M\in L^{\infty}(\Omega)$ and $\Lambda(M) := \sup  Re(\sigma(K
+MI)) >0$. 

Then  every nonegative nontrivial equilibrium is strictly positive and

\noindent i) For any nonnegative nontrivial $u_{0}\geq 0$ there exists
a positive  equilibria $\psi$ such that 
\begin{displaymath}
\liminf_{t\to \infty}  u(x, t, u_{0}) \geq \psi (x), \quad
a.e. \quad x\in \Omega .
\end{displaymath}
For such $\psi$ there exists some positive initial data $0<v_{0} <
\psi$ such that $\lim_{t\to \infty}  u(t,  v_{0}) = \psi$, 
in $L^{q}(\Omega)$ for any $1\leq q < \infty$.

\noindent ii) Moreover, if $\Lambda(M)$     is a principal eigenvalue  (see e.g.
(\ref{eq:measure_h=m}), (\ref{eq:1_over_h-m_notintegrable}), 
(\ref{eq:oscilation_h})),  there
exists an strictly positive   equilibrium $0< \varphi_{m}^{++} \leq 
\varphi_{M}$ such that it is minimal among the positive equilibria and
 stable from below as in Corollary \ref{extrem_equil_stable_from_above}.

 The assumption above holds in particular if  
\begin{displaymath}
  \lim_{s\to 0} \frac{f(x,s)}{s} = n(x), \quad \mbox{uniformly in
    $\Omega$ and} \quad 
  \Lambda (n) >0  \ \mbox{where}  \ n\in L^{\infty}(\Omega). 
\end{displaymath}

\end{proposition}
\begin{proof}
That every nonegative nontrivial equilibrium is strictly positive
comes from Proposition \ref{varphi_M_is_positive}. 


\noindent i) 
Fix $0< \tilde{\lambda} < \Lambda$ and, using Remark
\ref{rem:lower_bounds}, 
chose  $0<
\tilde{\varphi} \in L^{\infty}(\Omega)$  such that $\tilde{\lambda}
\tilde{\varphi} < K\tilde{\varphi} + M(x)\tilde{\varphi}$  
and $0< \tilde{\varphi} \leq  1$. 


Then observe that for  $\phi = \gamma \tilde{\varphi}$ with $0<\gamma \leq s_{0}$  
\begin{displaymath}
  K\phi + f(x,\phi) \geq  K\phi + M(x) \phi \geq \tilde{\lambda} \phi \geq 0.
\end{displaymath}
Hence $\phi$ is a subsolution of (\ref{nonlinear_loc_lip_asymp_est})
and then $  \phi \leq u(t, \phi)$ for $t\geq 0$, 
and $\phi \leq \Phi$ implies  that $u(t, \phi)$ is increasing and bounded by
$\Phi$, so it
converges  in $L^{q}(\Omega)$ for any $1\leq q < \infty$ to a bounded
positive equilibria. Denote this limit $u_{\gamma}$
with $0< \gamma \leq s_{0}$. Then  $\gamma \tilde{\varphi} \leq
u_{\gamma}$ and is uniformly bounded
in $L^{\infty}(\Omega)$ and is increasing in
$\gamma$.

If $u_{0} \geq 0$ is nontrivial we know from
Corollary \ref{nonlinear_max_pple_f_Loc_Lipschitz}
 that $\inf_{\Omega}u(t,u_{0})  >0$ for all
$t>0$. Hence we can assume without loss of generality that $u_{0} \geq
\alpha >0$ with $\alpha <
s_{0}$ and  then $u_{0} \geq \alpha \tilde{\varphi}$ which gives $
u(t, u_{0}) \geq u(t, \alpha \tilde{\varphi}) \geq  \alpha 
  \tilde{\varphi}$ for $t\geq 0$ 
and the term in the middle converges to the equilibria
$u_{\alpha}$. Hence
$\liminf_{t\to \infty}  u(x, t, u_{0}) \geq u_{\alpha} (x)$ 
a.e.  $x\in \Omega $. 
So $\psi = u_{\alpha}$ satisfies the statement i).


%
%

\noindent ii) 
%
Take $\tilde{\lambda} = \Lambda$ in part  i)  and $\tilde{\varphi} =
\varphi$  a positive bounded eigenfunction associated to the principal
eigenvalue, $\Lambda= \Lambda(M)$, of the operator $K + MI$, which 
is simple and moreover $\varphi \geq \alpha >0$ in $\Omega$ and
normalized $\|\varphi \|_{L^{\infty}(\Omega)} = 1$. 
Again, for  $0<\gamma \leq s_{0}$ and $\phi = \gamma \varphi \leq s_{0}$ then 
\begin{displaymath}
  K\phi + f(x,\phi) \geq K\phi + M(x) \phi = \Lambda \phi \geq 0 . 
\end{displaymath}
Hence $\phi$ is a subsolution of  (\ref{nonlinear_loc_lip_asymp_est}). 
Arguing as above we have that $u$  
converges  in $L^{q}(\Omega)$ for any $1\leq q < \infty$ to a bounded  positive equilibria. Denote this limit $u_{\gamma}$
with $0< \gamma \leq s_{0}$. Then  $u_{\gamma}$ is uniformly bounded
in $L^{\infty}(\Omega)$ and is increasing in
$\gamma$. Then the monotonic  limit $u_{*} = \lim_{\gamma \to 0} u_{\gamma}$
exists in $L^{q}(\Omega)$ for any $1\leq q < \infty$ and passing to the limit in $Ku_{\gamma} +
f(x,u_{\gamma})=0$ we obtain $Ku_{*} + f(x,u_{*}) =0$.
Below we show that $u_{*}$ is nontrivial (hence strictly positive), minimal and stable
from below.

Denote $F(x,s) =
  \begin{cases}
    M(x) s, & 0\leq s \leq s_{0} \\ f(x,s) & s> s_{0}
  \end{cases}$ 
so  $F(x,s) \leq f(x,s)$ for all $s \geq 0$ and globally  Lipschitz. Then if $0\leq u_{0}$ is nontrivial, we get 
\begin{displaymath}
  u_{f}(t,u_{0}) \geq u_{F}(t, u_{0})  \quad t\geq 0. 
\end{displaymath}
In particular if $u_{0} = \phi = \gamma \varphi$ with $0<\gamma < s_{0}$ we have, by
uniqueness, $u_{F}(t,u_{0}) = \gamma \varphi e^{\Lambda t}$ and is
increasing in time for as
long as $ \gamma  e^{\Lambda t} \leq s_{0}$, that is, for $t\leq
t_{0} (\gamma)= \frac{1}{\Lambda} \log(\frac{s_{0}}{\gamma})$. Also
$u_{F}(t,u_{0})$ is increasing in time since it is increasing for
$0\leq t\leq t_{0}$.
In particular, $u_{f}(t, \phi) \geq u_{F}(t,
u_{0})  \geq u_{F}(t_{0}, u_{0}) = s_{0} \varphi$, for  $t\geq
    t_{0}(\gamma)$. 


As above we can assume without loss of generality that we take initial data such that  $u_{0} \geq
\alpha >0$ and even more that we chose $\gamma< s_{0}$ such that $u_{0} \geq
\gamma \varphi$. Therefore $  u_{f}(t, u_{0}) \geq s_{0} \varphi$, for
$t\geq    t_{0}(\gamma)$.

In particular  every nontrivial   equilibrium $\psi$ 
satisfies $\psi \geq  s_{0} \varphi$ and then for every $\gamma$ we
have  $u_{\gamma} \geq  s_{0} \varphi$. Hence $u_{*} \geq
s_{0} \varphi$ and is nontrivial.


Also from    $\psi  \geq
\gamma \varphi$, we get $\psi=u_{f}(t, \psi) \geq u_{f} (t, \gamma
\varphi) \to u_{\gamma}$ as $t\to \infty$. Therefore $\psi \geq u_{*}$
and $u_{*}$ is the minimal positive equilibrium. 

To conclude  we show that $u_{*}$ is stable  from
below. In fact if $u_{*} \geq u_{0} >\alpha >0$, chose $0<\gamma <
s_{0}$ such that $u_{0} \geq \gamma \varphi$. Then for $t\geq 0$,
$u_{*} = u_{f}(t,u_{*}) \geq u_{f}(t,u_{0}) \geq u_{f}(t, \gamma
\varphi) $ 
and taking limit as $t\to \infty$ we know $u_{f}(t, \gamma
  \varphi) \to u_{\gamma} \geq u_{*}$ hence $u_{f}(t, u_{0}) \to
  u_{*}$.  

Finally, if $  \lim_{s\to 0} \frac{f(x,s)}{s} = n(x)$ uniformly in
$\Omega$ and $ \Lambda (n) >0$, 
since   $n \in L^{\infty}(\Omega)$ then  for $\eps>0$ small there
exists $s_{0}$ such that for $0\leq 
s\leq s_{0}$ and $x\in \Omega$, $f(x,s) \geq (n(x) - \eps) s$. 
Also,
for $\eps>0$ small, $  \Lambda ( n-\eps ) =  \Lambda
  ( n) -\eps >0$. 
Hence the assumptions are  satisfied with $M(x) = n(x)
-\eps$. Observe that $ \Lambda ( n-\eps )$ is a principal eigenvalue
iff  $ \Lambda( n)$ is so.  
\end{proof}

In particular we get the following result in the spirit of
Brezis--Oswald, \cite{brezis86:_remar}. 

\begin{proposition} 
    Assume $f:\Omega \times (0,\infty)$ is locally Lipschitz,
    $f(x,0)=0$, 
  \begin{displaymath}
    f(x,s) \leq C s + D, \qquad s\geq 0, \quad x\in \Omega
  \end{displaymath}
for some constants $C,D \geq 0$ and such that the limits 
\begin{displaymath}
 M_{0}(x) = \lim_{s\to 0^{+}} \frac{f(x,s)}{s}, \quad  M_{\infty}(x) = \lim_{s\to \infty} \frac{f(x,s)}{s}, \quad 
\end{displaymath}
exist uniformly in $\Omega$ and $M_{0}, M_{\infty} \in
L^{\infty}(\Omega)$ are such that $  \Lambda (M_{\infty}) < 0 <
\Lambda(M_{0})$. 

Then there exists at least a positive and bounded  solution of 
\begin{displaymath}
  Ku(x)\!+\!f(x,u(x))=0,  \quad x\in\Omega. 
\end{displaymath}

\end{proposition}
\begin{proof}
Assumption  $ \Lambda (M_{\infty}) < 0  $ implies that we have
(\ref{structure_condition_on_f}) 
with $\sup Re\Big(\sigma_{X}(K+ CI ) \Big)\leq  -\delta<0$
and we can apply Proposition \ref{asymp_estimate_phi_nonnegative} and
Theorem \ref{two_extremal_equilibria}. 

Assumption $0 < \Lambda(M_{0}) $ implies we can use Proposition
\ref{prop:minimal_positive_0_unstable}. 
\end{proof}

The next result gives a sufficient condition for uniqueness of
positive equilibria.


\begin{theorem}
\label{thr:global_asymptot_stable}
  
As in  Proposition \ref{varphi_M_is_positive},
assume $g(x)= f(x,0)\ge 0$,  $x\in \Omega$ and  additionally that
the kernel $J$ is symmetric, that is, $J(x,y) = J(y,x)$ and 
	\begin{equation}\label{f/s_decreasing}
		\displaystyle \frac{f(x,s)}{s}\;\mbox{ is
                  decreasing  in $s$ for a.e. $x\in \Omega$. }  
	\end{equation}
	
	Then $\varphi_M$ is the  unique nontrivial nonnegative
        equilibrium of (\ref{nonlinear_loc_lip_asymp_est}),  and for
        every nontrivial $u_{0}\geq 0$ we have
        \begin{displaymath}
          \lim_{t \to \infty} u(t,u_{0}) = \varphi_{M}. 
        \end{displaymath}
That is, $\varphi_{M}$ 
        is  globally asymptoticaly for the solutions
        of (\ref{nonlinear_loc_lip_asymp_est})  with nonnegative initial data.

Moreover 

\noindent i) If   $g$ is not identically zero then $\varphi_{M}$ is
strictly positive in $\Omega$.

\noindent ii) If  $g =0$ assume 
\begin{displaymath}
  \lim_{s\to 0} \frac{f(x,s)}{s} =  n(x), \quad \mbox{uniformly in
    $\Omega$, with $n\in L^{\infty}(\Omega)$.} 
\end{displaymath}

Then, if  $  \Lambda ( n ) \leq 0$ 
we have  $\varphi_{M} =0$, while if $  \Lambda ( n  ) > 0$
then $\varphi_{M}$ is strictly positive in $\Omega$.

\end{theorem}
\begin{proof}
  	From Theorem \ref{two_extremal_equilibria}, let $\varphi_M\in 
	L^{\infty}	(\Omega)$ 
	be the maximal equilibria  of
        (\ref{nonlinear_loc_lip_asymp_est}) and assume $\varphi_{M}$
        is nontrivial (otherwise there is nothing to prove). Now, assume that 
	$\psi$ is another nontrivial nonnegative equilibria,  
	then $0\leq \psi\le \varphi_M$. Thus, $\psi\in
        L^{\infty}(\Omega)$ and from Proposition 
	\ref{varphi_M_is_positive}, 	$0<\psi \leq \varphi_{M}$. 


Notice that the equation for any positive equilibrium   $\xi$ can be
written as $K\xi  +
\frac{f(\cdot, \xi)}{\xi} \xi=0$. Hence, since  $\psi, \varphi_{M}$
are positive equilibria we get that 
\begin{displaymath}
  \Lambda \left( \frac{f(\cdot, \psi)}{\psi}\right) = 0 =   \Lambda \left( \frac{f(\cdot, \varphi_{M})}{\varphi_{M}}\right)
\end{displaymath}
and they are principal eigenvalues (with  $\psi, \varphi_{M}$ as
positive eigenfunctions respectively)  while, at the same time, since
$0<\psi \leq \varphi_{M}$, 
\begin{displaymath}
  \frac{f(\cdot, \psi)}{\psi} \geq  \frac{f(\cdot, \varphi_{M})}{\varphi_{M}}
\end{displaymath}
with strict inequality in a set of positive measure. Now multiplying
$K(\varphi_{M}) + \frac{f(\varphi_{M})}{\varphi_{M}} \varphi_{M}=0$ by $\psi$ and  $K(\psi) + \frac{f(\psi)}{\psi}\psi
=0$ by $\varphi_{M}$,  integrating in $\Omega$ and using that $J$ is 
symmetric, we have $\int_{\Omega} K(\varphi_{M}) \psi = \int_{\Omega}
\varphi_{M} K (\psi)$ and then 
\begin{displaymath}
  \int_{\Omega} \Big( \frac{f(\psi)}{\psi}-
  \frac{f(\varphi_{M})}{\varphi_{M}}\Big)  \varphi_{M}\psi =0 
\end{displaymath}
and this is a contradiction since the integrand is nonnegative and
nonzero in a set of positive measure. 
Therefore $\varphi_{M}$ is the unique nonnegative equilibrium.

If   $g$ is not identically zero, then by Proposition
\ref{varphi_M_is_positive} we get that $\varphi_{M}$ is strictly
positive in $\Omega$. 

On the other hand, assume  $g=0$.  If  $\Lambda (n) >0$ by
Proposition \ref{prop:minimal_positive_0_unstable} we also get that $\varphi_{M}$ is strictly
positive in $\Omega$. 

If $\Lambda (n)< 0$ by (\ref{f/s_decreasing}) we get
$\frac{f(x,s)}{s} \leq  \lim_{s\to 0} \frac{f(x,s)}{s} =  n(x)$, 
 $s>0$, $x\in \Omega$ 
and then $f(x,s) \leq n(x) s$ for $s>0$, $x\in \Omega$. 
That is, $f$ satisfies (\ref{structure_condition_on_f}) for $s>0$, $C(x)=
n(x)$ and $D(x)=0$. Hence (\ref{semigroup_LC_decays}) is satisfied and
in (\ref{equation_PHI}) we get $\Phi=0$ and then
$\varphi_{M}=0$. 

Finally, if $\Lambda(n)=0$ and $\varphi_{M}$ was nontrivial and thus
strictly positive, we will get as above $\Lambda \left( \frac{f(\cdot,
  \varphi_{M})}{\varphi_{M}}\right) =0$ and is a principal eigenvalue and at
the same time $\frac{f(\cdot,  \varphi_{M})}{\varphi_{M}} \leq n$
which again contradicts the strict
monotonicity of the principal eigenvalue in Proposition 2.36 in
\cite{RB2017_max_pples}. 
\end{proof}

\begin{remark}
%
%
\noindent i)
Theorem \ref{thr:global_asymptot_stable} is known to hold for local
diffusion problems like (\ref{eq:local_RD}), see
e.g. \cite{Anibal_Alejandro_extremal_equilibria} and references
therein, assuming in (\ref{f/s_decreasing}) only that
$\frac{f(x,s)}{s}$ is nonincreasing. The reason for this is that for
(\ref{eq:local_RD}) the maximum principle   implies that   in the
proof of the theorem, $\{ \psi <\varphi_{M}\} = \Omega$ while in 
the case of nonlocal problems in this paper we can not guarantee this,
see Example \ref{example_discontinuous_equilibria}. In such a case,
without strict decreasing in (\ref{f/s_decreasing}) we can not
conclude in the proof above that $  \frac{f(\cdot, \psi)}{\psi} \geq
\frac{f(\cdot, \varphi_{M})}{\varphi_{M}}$ 
with strict inequality in a set of positive measure, see
(\ref{eq:4_uniqueness}) for the case of logistic nonlinearities
below. 

\noindent ii) 
  If $f(x,s)$ is regular in $s$, we have that (\ref{f/s_decreasing}) 
  is equivalent to  $f(x,s) > s \frac{\partial f}{\partial s}(x,s)$,  $s\geq 0$, $x\in \Omega$. 
 This holds in
particular if  $f(x,s)$ is strictly concave in $s$ since by the mean value
theorem,  for some $0\leq \xi(x) \leq s$ we have $\frac{f(x, s)}{s} \geq     \frac{f(x, s)-f(x,0)}{s} =
   \frac{\partial f}{\partial s}(x,\xi) >  \frac{\partial f}{\partial
     s}(x,s)$. 

\end{remark}




\subsection{Logistic type nonlinearities}
\label{logistic_nonlin}

Consider  logistic nonlinearities 
\begin{displaymath}
  f(x,s) = g(x) + n(x) s - m(x) |s|^{\rho-1} s
\end{displaymath}
with $g, n, m\in L^{\infty}(\Omega)$, $m\ge 0$ not identically zero
and $\rho>1$.

Then we have
\begin{equation} \label{eq:first_bound_logistic_f}
  f(x,s)s \leq n(x) s^{2} + |g(x)||s| 
\end{equation}
and, since  $g, n \in L^{\infty}(\Omega)$, from Proposition \ref{global_sol_initial_data_bounded_2} we have
existence and uniqueness solution of
\eqref{nonlinear_loc_lip_asymp_est} for $u_{0}\in
X=L^{\infty}(\Omega)$ or $X=C_b(\Omega)$.

Moreover since 
\begin{displaymath}
  \frac{\partial f}{\partial s}(x,s) = n(x) -\rho m(x) |s|^{\rho-1}
  \leq  n(x),  \quad
  \left|\frac{\partial f}{\partial s}(x,s)\right| \leq  c(1+ |s|^{\rho-1}).
\end{displaymath}
Then, if  $|\Omega|<\infty$, $J(x,y)=J(y,x)$,  from Theorem
\ref{global_solutions_f_growth} and we have existence and uniqueness
of solution of \eqref{nonlinear_loc_lip_asymp_est} for $u_{0}\in
L^{p\rho}(\Omega)$, $1\leq p <\infty$.  

As for  asymptotic estimates and extremal equilibria results we have
the following results. First, in the following result the asymptotic
behavior of solutions is determined by the linear terms in the
equation.

\begin{proposition}

  Assume  $  \Lambda(n) <0$.

  \noindent i)
  Then there exist  two bounded extremal equilibria, $\varphi_m\le 	\varphi_M $,  which enclose the
asymptotic behavior of the solutions and that are stable from  below and
from above.

\noindent ii)
For nonnegative solutions,    if  $g \geq 0$ and $g\neq 0$,  there exists a 
minimal positive  equilibria which is 
stable from below  for nonnegative solutions. On the other hand, if
$g=0$ then all nonnegative solutions
converge to zero.

If 
\begin{equation}\label{eq:4_uniqueness}
  |\{g=0\} \cap \{m=0\}| =0 . 
\end{equation}
and $g\neq 0$, then there exists a unique positive equilibria which is
moreover  globally asymptotically
for nonnegative initial data. 

\end{proposition}
\begin{proof}
  
From
(\ref{eq:first_bound_logistic_f}) we can take $C(x)=n(x)$, $D(x)=|g(x)|$.
Hence, from  Theorem \ref{two_extremal_equilibria} and
  Corollary \ref{extrem_equil_stable_from_above} we  
obtain the existence of two bounded extremal equilibria, $\varphi_m\le 	\varphi_M $,  which enclose the
asymptotic behavior of the solutions and that are stable from  below and from above.

For nonnegative solutions, if  $g \geq 0$ and $g\neq 0$,  from Proposition
\ref{varphi_M_is_positive} we have existence of a
minimal positive  equilibria which is 
stable from below  for nonnegative solutions. On the other hand, if
$g=0$ then $\varphi_{M}=0= \Phi$ and all nonnegative solutions
converge to zero.

Notice we  also have  
\begin{displaymath}
  {f(x,s)\over s} = {g(x)\over s}+n(x)-m(x)|s|^{\rho-1} 
\end{displaymath}
and then (\ref{f/s_decreasing}) holds provided (\ref{eq:4_uniqueness})
holds true. 
In such a case,  from Theorem
\ref{thr:global_asymptot_stable}, if $g\neq 0$, then $\varphi_{M}>0$
is the unique nonnegative equilibria and is globally asymptotically
for nonnegative initial data. 
%
\end{proof}

\begin{proposition}
  
Assume  $  \Lambda(n) \geq 0$.

\noindent i) Then there exist  two bounded extremal equilibria, $\varphi_m\le 	\varphi_M $,  which enclose the
asymptotic behavior of the solutions and that are stable from  below and
from above.

\noindent ii) For nonnegative solutions,    if  $g \geq 0$ and $g\neq 0$,  there exists a 
minimal positive  equilibria which is 
stable from below  for nonnegative solutions. 

Additionally, 

ii.1) Assume $m(x) \geq m_{0} >0$ in $\Omega$.



If $g\neq 0$ or $g=0$ and  $\Lambda(n)>0$, then 
there exists a unique positive equilibria which is
moreover  globally asymptotically
for nonnegative initial data. 
Finally, if $g=0$ and $\Lambda(n)=0$
then every nonnegative solution converges to
zero.

 ii.2)
Assume  $m(x)$ vanishes in a set of
positive measure of $\Omega$.


If $g=0$ then $\varphi_{m} \leq 0 \leq \varphi_{M}$.  
Moreover, if  $g\neq 0$ is such that    (\ref{eq:4_uniqueness}) holds,
then  there exists a  unique nonnegative equilibria and is globally asymptotically
for nonnegative initial data.

\end{proposition}
\begin{proof}
\noindent i)
Assume $m(x) \geq m_{0} >0$ in $\Omega$. 
Now choose $A >0$ such that $  \Lambda(n-A) = \Lambda(n) -A  <0$  
and write 
\begin{displaymath}
  f(x,s)s \leq  |g(x)||s| + \big(n(x)-A\big)  s^{2} +  |s|\big( A|s| -
  m(x)|s|^{\rho}\big)  . 
\end{displaymath}
Then, Young's inequality gives, for any $\eps>0$, 
\begin{displaymath}
  A|s| -  m(x)|s|^{\rho} \leq A|s| -  m_{0}|s|^{\rho} \leq \eps
  |s|^{\rho} + C_{\eps} A^{\rho'} -    m_{0}|s|^{\rho}  
\end{displaymath}
and taking $\eps= \frac{m_{0}}{2}$ we get
\begin{displaymath}
  f(x,s) s \leq \big(n(x)-A\big)  s^{2} +  \big( |g(x)| + C_{\eps} A^{\rho'} \big) |s| .
\end{displaymath}
Then  we can take $C(x)= n(x)-A$ and
$D(x)= |g(x)| + C_{\eps} A^{\rho'}  \in L^{\infty}(\Omega)$ and from
Theorem \ref{two_extremal_equilibria} and Corollary  
\ref{extrem_equil_stable_from_above} we have again  the existence of
two bounded extremal equilibria, $\varphi_m\le 	\varphi_M $.

If $g\ge 0$ and $g\neq 0$, we obtain
from Proposition \ref{varphi_M_is_positive} the existence of a minimal
positive equilibria $0< \varphi_{m}^{+} \leq \varphi_{m}$, stable from
below for nonegative equilibria. Also,  if
$g=0$ then $\varphi_{m} \leq 0 \leq \varphi_{M}$.  

Observe that if $g=0$ then 
\begin{displaymath}
  \lim_{s\to 0} \frac{f(x,s)}{s} =  n(x), \quad \mbox{uniformly in
    $\Omega$}. 
\end{displaymath}
Therefore since   (\ref{eq:4_uniqueness}) holds, then  (\ref{f/s_decreasing})
holds and from   Theorem \ref{thr:global_asymptot_stable}, 
if $g\neq 0$ or $g=0$ and  $\Lambda(n)>0$, then $\varphi_{M}>0$
is the unique nonnegative equilibria and is globally asymptotically
for nonnegative initial data.
Finally, if $g=0$ and $\Lambda(n)=0$
then $\varphi_{M}=0$ and every nonnegative solution converges to
zero.

\noindent ii) 
Now we consider the case in which $m(x)$ vanishes in a set of
positive measure of $\Omega$. 
Assume that for some $\delta >0$,  $\omega'
\subset \supp(m)$ is such that $m(x)>\delta$ for $x\in
\omega'$ and of positive measure and  denote  $\Omega'= \Omega \setminus
\omega'$. 

Now if $x\in\Omega'$ then from (\ref{eq:first_bound_logistic_f})  we take $C(x)=n(x)$ and $D(x)=|
g(x)|$ for  $x\in\Omega'$.
On the other hand, if  $x\in \omega'$,  we proceed as in  part i)  above, with a large  $A>0$ 
\begin{displaymath}
  f(x,s)s \leq \big(n(x)-A\big)  s^{2} +  \big( |g(x)| + C_{\eps}
  A^{\rho'} \big) |s|, \qquad x\in \omega' . 
\end{displaymath}

Then we have (\ref{structure_condition_on_f}), 
that is, $f(x,s)s \le C(x)s^2+D(x)|s|,\qquad  s \in \R, \quad x \in
        \Omega $ 
with 
\begin{displaymath}
  D(x) =
\begin{cases}
  |g(x)|, & x\in \Omega' \\
|g(x)|+C_{\eps}A^{\rho'},  & x\in \omega'
\end{cases}, 
\quad 
C(x) =
\begin{cases}
  n(x), & x\in \Omega' \\
n(x)-A, & x\in \omega'
\end{cases}
\end{displaymath}
and by Proposition \ref{prop:shifting_the_spectrum} below, with   $A$ large enough 
we have  $\Lambda(C)<0$. 

Hence, again from
Theorem \ref{two_extremal_equilibria} and Corollary  
\ref{extrem_equil_stable_from_above} we have again  the existence of
two bounded extremal equilibria, $\varphi_m\le 	\varphi_M $. If $g\ge 0$ and $g\neq 0$, we obtain
from Proposition \ref{varphi_M_is_positive} the existence of a minimal
positive equilibria $0< \varphi_{m}^{+} \leq \varphi_{m}$, stable from
below for nonegative equilibria. Also,  if
$g=0$ then $\varphi_{m} \leq 0 \leq \varphi_{M}$.

Moreover, if  $g\neq 0$ is such that    (\ref{eq:4_uniqueness}) holds, then  (\ref{f/s_decreasing})
holds and from   Theorem \ref{thr:global_asymptot_stable},  then $\varphi_{M}>0$
is the unique nonnegative equilibria and is globally asymptotically
for nonnegative initial data.
\end{proof}

Now we prove the result used above that  states that by acting on an arbitrary small
subset of the domain with a large negative constant,  we can shift the spectrum of a nonlocal operator
$K+hI$ to have negative real part. Observe that this result is not
true for local diffusion operators $-\Delta + hI$, see
\cite{Anibal_Alejandro_extremal_equilibria}.

First, define, if  $\Omega ' \subset \Omega$ and for $\varphi \in X(\Omega')=
 L^{p}(\Omega')$, with $1\leq p\leq \infty$,  the nonlocal
 operator in $X(\Omega')$ 
\begin{displaymath}
  K_{\Omega'} \varphi (x) = \int_{\Omega'} J(x,y) \varphi(y) \,
  dy, \quad x \in \Omega ' .
\end{displaymath}
and for $h\in L^{\infty}(\Omega')$ denote $\Lambda(h, \Omega') = \sup
Re(\sigma(K _{\Omega'} +h I)$. 
Also, denote $\omega ' =
\Omega\setminus \Omega '$.

\begin{proposition}
  \label{prop:shifting_the_spectrum}

  Assume $|\omega'| >0$ and $h \in L^{\infty}(\Omega)$. 
 Then for sufficiently large $A>0$ and defining 
$H(x)=
\left\{\begin{array}{ll}
         h(x) &  x\in\Omega',
         \\
h(x) - A &x\in\omega', 
\end{array}
\right.
$ 
 for the operator $K+HI$ in $X=L^{p}(\Omega)$, $1\leq p\leq \infty$, we have $\Lambda (H)
<0$. 
\end{proposition}
\begin{proof}
For $\varphi \in X= L^{p}(\Omega)$ we have,  
if $x\in \Omega'$ 
\begin{displaymath}
   \Big({K\varphi + H\varphi\over \varphi} \Big) (x) =
   \frac{K(\varphi  \chi_{\omega'}) (x)  +  K_{\Omega'}(\varphi 
  \chi_{\Omega'})(x) +h(x) \varphi (x)}{\varphi(x)}
\end{displaymath}
hence
\begin{displaymath}
  \sup_{\Omega'} {K\varphi + H\varphi\over \varphi} \leq
  \sup_{\Omega'} \frac{K(\varphi  \chi_{\omega'})} {\varphi} +  \sup_{\Omega'}  \frac{K_{\Omega'}(\varphi  
  \chi_{\Omega'}) + h \varphi}{\varphi} .
\end{displaymath}

If $x\in \omega'$ 
\begin{displaymath}
   \Big({K\varphi +H\varphi\over \varphi} \Big) (x) =
   \frac{K(\varphi  \chi_{\Omega'}) (x)  +  K_{\omega'}(\varphi 
  \chi_{\omega'})(x) + (h(x) -A) \varphi (x)  \chi_{\omega'}}{\varphi(x)}
\end{displaymath}
hence
\begin{displaymath}
  \sup_{\omega'} {K\varphi + H\varphi\over \varphi} \leq
  \sup_{\omega'} \frac{K(\varphi  \chi_{\Omega'})} {\varphi} +  \sup_{\omega'}  \frac{K_{\omega'}(\varphi  
  \chi_{\omega'}) + (h -A) \varphi}{\varphi} .
\end{displaymath}

Therefore, using Theorem \ref{thr:Lambda}, 
\begin{displaymath}
  \Lambda (H)=\inf\limits_{0<\varphi \in X}\sup\limits_{\Omega}
  {K\varphi + H \varphi\over \varphi} \leq  \inf\limits_{0<\varphi \in
    X} \Big(  \sup_{\Omega'} \frac{K(\varphi  \chi_{\omega'})}
  {\varphi}  +  \sup_{\omega'} \frac{K(\varphi  \chi_{\Omega'})}
  {\varphi} \Big)  + \Lambda(h, \Omega') + \Lambda(h-A, \omega') . 
\end{displaymath}

Now, taking $\varphi =1$,
\begin{displaymath}
  \inf\limits_{0<\varphi \in
    X} \Big(  \sup_{\Omega'} \frac{K(\varphi  \chi_{\omega'})}
  {\varphi}  +  \sup_{\omega'} \frac{K(\varphi  \chi_{\Omega'})}
  {\varphi} \Big)   \leq \sup_{\Omega} K(1) = \sup_{\Omega} h_{0}. 
\end{displaymath}
and since
$\Lambda(h-A, \omega') =
\Lambda(h, \omega') -A$ we get 
\begin{displaymath}
   \Lambda (H) \leq \sup_{\Omega} h_{0} + \Lambda(h, \Omega') + \Lambda(h, \omega') -A <0 
\end{displaymath}
for sufficiently large $A$.  
\end{proof}


\section{Further comments on compactness and asymptotic behavior}
\label{asymp_beh}

Observe that from the asymptotic estimates in Section
\ref{asymptotic_estimates},  in order to study the asymptotic dynamics
of (\ref{intro_nonlinear_diffusion}) it is enough to take (bounded)  initial data in
the set   
\begin{displaymath}
  \mathcal{B}= \{u_{0}, \ |u_{0}(x)|\leq \Phi(x), \ x\in \Omega\}
  \subset L^{\infty}(\Omega) 
\end{displaymath}
 and assume thereafter that $f(x,s)$ is globally  Lipschitz in the second
variable. Notice that from Proposition
\ref{asymp_estimate_phi_nonnegative} this set of initial data is
invariant for (\ref{intro_nonlinear_diffusion}), that is, the nonlinear semigroup
(\ref{eq:nonlinear_semigroup}) satisfies 
\begin{displaymath}
  S(t) \mathcal{B} \subset \mathcal{B}, \ t \geq 0. 
\end{displaymath}

For such class of initial data, we
have that the semigroup $S(t)$ is  continuous   in the norm of 
$L^{p}(\Omega)$ for any  $1\leq p\leq \infty$, as in Section
\ref{exist_uniq_pos_Glob_Lipschitz}.  Also,  for $u_{0} \in \mathcal{B}$ we have that $u(t,u_{0}),
u_{t}(t,u_{0}), f(u(t,u_{0}))$ are uniformly bounded in
$L^{\infty}(\Omega)$,  for all $t\geq0$
and independent of $u_{0}$. In particular $S(\cdot)\mathcal{B}$ is
equicontinuous in $\mathcal{C}([0,\infty), L^{\infty}(\Omega))$. 

Notice that after the results in Section
\ref{section_extr_equil} we could as well reduce ourselves to take
initial data such that $\varphi_{m} \leq u_{0} \leq \varphi_{M}$ since
this order interval of bounded functions is also invariant for
(\ref{intro_nonlinear_diffusion}).

However, in contrast with  the standard diffusion equation
(\ref{eq:local_RD}),   (\ref{intro_nonlinear_diffusion}), has very
poor regularizing properties which 
makes the asymptotic behavior of solutions difficult to define and
analyze. For example,  for linear problems it was proved in 
Theorem 4.5 in \cite{ARB_SSG_16} that if $h(x) \geq \alpha >0$ then
$e^{(K-hI)t}$ is asymptotically smooth. This weak compactness does not seem
enough to translate any compactness to the nonlinear semigroup
(\ref{eq:nonlinear_semigroup}) given by the variations of constants
formula 
since the Nemitcky operator $f$ does not have any compactness
properties between the Lebesgue spaces. Hence the semigroup 
\begin{displaymath}
   S(t): \mathcal{B} \to \mathcal{B} , \ t\geq 0
\end{displaymath}
is continuous but we lack of results to prove it is  compact, or asymptotically
compact. This precludes from having well defined $\omega$--limit sets
or an attractor describing the asymptotic behavior of solutions.

Also notice that if  $J(x,y) = J(y,x)$ then 
\begin{displaymath}
	  E(\varphi)=\frac{1}{2}\int_{\Omega}\int_{\Omega}J
	(x,y)(\varphi(y)-\varphi(x))^2 \, dy\,dx - \int_{\Omega} h_{0}(x)
	 \varphi^{2}(x)\, dx - \int_{\Omega} F(x,\varphi(x))\,
        dx 
	\end{displaymath}
where $F(x,s) = \int_{0}^{s} f(x,r)\, dr$ can also be assumed to be
globally bounded and Lipschitz,    is decreasing along
trajectories, that is 
\begin{displaymath}
  \frac{d}{dt} E(u(t,u_{0})) = -\int_{\Omega} |u_{t}(t,u_{0})|^{2} \leq 0
\end{displaymath}
and so it defines a strict  Lyapunov functional for the nonlinear
semigroup above.

If we had enough  compactness to guarantee that  for some $u_{0}$ as
above and for some sequence 
$t_{n}\to \infty$, we have that $u(t_{n}, u_{0}) \to \xi$ a.e. in
$\Omega$, then we would have convergence  in $L^{p}(\Omega)$ for  $1\leq p<\infty$  
since all functions involved are in $\mathcal{B}$. Then  $\xi$ is necessarily
an equilibria.  
However the lack of compactness/smoothing  mentioned above
precludes from guaranteeing that the trajectory $u(\cdot,u_{0})$
accumulates somewhere a.e. in $\Omega$,  as $t\to \infty$.

 We could also  consider the set $\mathcal{B}$ endowed with the
weak convergence in, say, $L^{2}(\Omega)$, which we denote
$\mathcal{B}_{w}$, which is a closed, convex,  compact (hence
complete) metric space. From the bounds  above on the semigroup  and
Ascoli--Arzela's theorem  we get
that $S(\cdot)\mathcal{B}_{w}$ is relatively compact in  
$\mathcal{C}([0,\infty), \mathcal{B}_{w})$. Again  if for a weakly convergent
sequence $u_{0}^{n} \to u_{0}$ we had $u(\cdot,u_{0}^{n}) \to \xi$
a.e. $(0,T) \times \Omega$, for each $T>0$,  then we would get that $\xi = u(\cdot,
u_{0})$ and the semigroup 
\begin{displaymath}
  S(t): \mathcal{B}_{w} \to \mathcal{B}_{w} , \ t\geq 0
\end{displaymath}
would be continuous and, obviously, compact. In such a case, the semigroup would have a global
attractor in $\mathcal{B}_{w}$. However, again the poorly regularizing effect
of the nonlocal diffusion equations seem not enough to prove the
pointwise a.e. convergence required in the argument above.

\addcontentsline{toc}{chapter}{Bibliography}

\end{document}